\documentclass[12pt]{amsart}

\usepackage{enumerate,amssymb,amsmath,mathrsfs}
\usepackage{amssymb,amsfonts,amsthm,amscd}
\usepackage[colorlinks=true,linkcolor=blue,citecolor=blue]{hyperref}
\usepackage{cite}
\usepackage[mathscr]{eucal}     
\usepackage{graphicx}
\usepackage[arrow, matrix, curve]{xy}
\usepackage{color}
\usepackage[margin=3cm]{geometry}
\definecolor{gray}{gray}{.75}
\definecolor{gray2}{gray}{.50}

\usepackage{tikz}

\newcommand{\dom}{\mathcal{D}}

\newcommand{\ccl}{\textup{ccl}}

\newcommand{\Z}{\mathbb{Z}}
\newcommand{\N}{\mathbb{N}}
\newcommand{\R}{\mathbb{R}}
\newcommand{\C}{\mathbb{C}}

\newcommand{\A}{\mathrm{\alpha}}
\newcommand{\w}{\mathrm{\omega}}

\newcommand{\e}{\epsilon}

\newtheorem{thm}{Theorem}[section]
\newtheorem{cor}[thm]{Corollary}
\newtheorem{remark}[thm]{Remark}
\newtheorem{prop}[thm]{Proposition}
\newtheorem{lemma}[thm]{Lemma}
\newtheorem{defn}[thm]{Definition}

\numberwithin{equation}{section}

\begin{document}

\title[The metric anomaly of analytic torsion at the boundary of an even dimensional cone]
{The metric anomaly of analytic torsion at the boundary of an even dimensional cone}

\author{Boris Vertman}
\address{Mathematisches Institut, 
Universit\"at Bonn, 53115 Bonn, Germany}
\email{vertman@math.uni-bonn.de}
\urladdr{www.math.uni-bonn.de/people/vertman}

\thanks{The author was partially supported by the Hausdorff Center for Mathematics}


\begin{abstract}
{The formula for analytic torsion of a cone in even dimensions is comprised of three terms. 
The first two terms are well understood and given by an algebraic combination of the Betti numbers and 
the analytic torsion of the cone base. The third "singular" contribution is an intricate 
spectral invariant of the cone base. We identify the third term as the metric anomaly 
of the analytic torsion coming from the non-product structure of the cone at its regular boundary. Hereby 
we filter out the actual contribution of the conical singularity and identify the analytic torsion 
of an even-dimensional cone purely in terms of the Betti numbers and the analytic torsion of the cone base.}
\end{abstract}

\maketitle
\tableofcontents

\pagestyle{myheadings}
\markboth{\textsc{Analytic Torsion}}{\textsc{Boris Vertman}}

\section{Introduction and Statement of the Main Result}

Analytic torsion has been introduced by Ray and Singer in \cite{RaySin:RTA} 
as the analytic counterpart to the combinatorially defined Reidemeister-Franz torsion. 
The latter was the first topological invariant which was not a homotopy invariant, defined 
and studied by Reidemeister, Franz and de Rham in \cite{Re1}, 
\cite{Re2}, \cite{Fr} and \cite{Rh}. The equality between the analytic and Reidemeister-Franz torsion 
has been conjectured by Ray and Singer in \cite{RaySin:RTA}, and proved independently by Cheeger 
\cite{Che:ATA} and M\"uller \cite{Mu} for closed manifolds with a unitary representation of the 
fundamental group. 

In this article we study analytic torsion of an even-dimensional cone 
$\mathscr{C}(N)=(0,1]\times N$ over a closed odd-dimensional 
Riemannian manifold $(N,g^N)$, equipped with a warped product metric 
$g=dx^2\oplus x^2g^N$. The presented analysis provides the even-dimensional 
analogue of our results jointly with Werner M\"uller in \cite{MV}.

In the previous publication \cite{Ver:ATO}, we have 
derived an expression for the analytic torsion on $(\mathscr{C}(N),g)$
in terms of spectral invariants of $N$. The present article is devoted to an 
identification of one of these invariants as the metric anomaly of Br\"uning-Ma 
\cite{BruMa:AAF} for analytic torsion of $(\mathscr{C}(N),g)$ at the regular boundary. 

Overall we identify the analytic torsion of a cone in even dimensions it terms 
of cohomology and Ray-Singer torsion of the cone base, plus the metric anomaly 
at the regular boundary. This can be viewed as a step forward towards a Cheeger-M\"uller 
type result in the singular setup, compare the recent result by Hartmann-Spreafico \cite{HS2}.

\section*{Acknowledgements}

The author would like to express deep gratitude to Jeff Cheeger and 
Werner M\"uller for encouragement and many useful discussions that have led to the presented results. 
The author would also like to thank Stanford Department of Mathematics, where the results have been obtained, 
for hospitality and support. The author gratefully acknowledges the financial support of the German Research Foundation DFG, as well as the Hausdorff Institute at the University of Bonn.

\subsection{The Setup and Definition of Analytic Torsion}\label{section-torsion} 
Let $(N^n,g^N), n=\dim N,$ be an odd-dimensional smooth closed Riemannian manifold. 
Consider a bounded cone of even dimension over $N$
$$\mathscr{C}(N)=(0,1] \times N, \quad g=dx^2\oplus x^2g^N.$$ 
Consider a complex flat vector bundle $(E_N,\nabla_N,h_N)$ over $N$, induced by a unitary 
representation $\rho_N: \pi_1(N)\to U(n,\C)$ of the fundamental group of $N$, 
acting by deck transformations on the universal cover $\widetilde{N}$. 
As explicated in \cite{MV}, every flat complex vector bundle $(E,\nabla, h^E)$ 
over the cone $(\mathscr{C}(N),g)$ arises from such a flat complex vector bundle over its base $N$,
with $E=E_N\times (0,1]$ and $h^E \restriction \{x\}\times N=h_N$.
The flat covariant derivatives $\nabla$ and $\nabla_N$ 
are related for any section $s\in C^{\infty}([0,1], C^{\infty}(E_N))\equiv 
\Gamma(E)$ as follows
\begin{align}
\nabla s= \frac{\partial s}{\partial x}\otimes dx + \nabla_N s.
\end{align}
Let $(\Omega^*_0(\mathscr{C}(N),E),\nabla_*)$ be the associated twisted de Rham complex, 
where $\Omega^k_0(\mathscr{C}(N),E)$ are the $E$-valued differential $k$-forms with compact support in 
$\mathscr{C}(N)$, and $\nabla_k$ denotes the covariant differential induced by the flat connection $\nabla$ 
on $E$, densely defined with the domain $\Omega^k_0(\mathscr{C}(N),E)$. 
Let $\nabla_{k,\min}$ denote the graph closure of $\nabla_k$  and $\nabla_{k,\max}$ 
its maximal closed extension in $L^2(\Omega_0^*(\mathscr{C}(N),E),g,h^E)$.  
The relative and the absolute self-adjoint extensions of the Laplacian are defined by
\begin{equation}
\begin{split}
&\Delta_k^{\textup{rel}}:=\nabla^*_{k,\min}\nabla_{k,\min}+\nabla_{k-1,\min}
\nabla^*_{k-1,\min}, \\
&\Delta_k^{\textup{abs}}:=\nabla^*_{k,\max}\nabla_{k,\max}+\nabla_{k-1,\max}
\nabla^*_{k-1,\max}. 
\end{split}
\end{equation}

\begin{thm}[Cheeger \cite{Che:SGOlv:AAS}]
\label{heat-trace}
Let $(\mathscr{C}(N)=(0,1]\times N, g=dx^2\oplus x^2g^N)$ be an even-dimensional bounded cone over 
a closed oriented Riemannian manifold $(N^n,g^N)$. Let $(E,\nabla, h^E)$ be a flat complex 
Hermitian vector bundle over $\mathscr{C}(N)$. Then for the relative or absolute self-adjoint extension 
$\Delta_k^{\textup{rel/abs}}$ of the Laplacian on $\Omega^k_0(\mathscr{C}(N),E)$, we find
\begin{align}
\textup{Tr}\, e^{-t\Delta_k^{\textup{rel/abs}}} \sim \sum_{j=0}^{\infty}A_j t^{\frac{j-m}{2}}
+\sum_{j=0}^{\infty}C_j t^{\frac{j}{2}}+\sum_{j=0}^{\infty}G_j 
t^{\frac{j}{2}}\log t, \ t\to 0.
\end{align}
\end{thm}
The heat-trace expansion yields a meromorphic extension of $\Delta_k^{\textup{rel/abs}}$
\begin{align}
\zeta (s, \Delta_k^{\textup{rel/abs}}):=
\frac{1}{\Gamma(s)}\int_0^{\infty}t^{s-1}\textup{Tr}(e^{-t\Delta_k^{\textup{rel/abs}}}-
\mathscr{P}^{\textup{rel/abs}}_k) dt, \ \textup{Re}(s) \gg 0
\end{align}
to the whole complex plane, where $\mathscr{P}^{\textup{rel/abs}}_k$ denotes the orthogonal projection 
of $L^2\Omega^k_0(\mathscr{C}(N),E)$ onto the subspace $\mathscr{H}^k_{\textup{rel/abs}}(\mathscr{C}(N),E)$ 
of harmonic forms for $\Delta_k^{\textup{rel/abs}}$. 
The heat-trace expansion determines the behaviour of the zeta-function near $s=0$, and 
the following result is a direct consequence of Theorem \ref{heat-trace}. 

\begin{cor}\label{zeta}
\begin{align}
\Gamma(s)\zeta (s, \Delta_k^{\textup{rel/abs}}) \sim \sum_{j=0}^{\infty}A_j\frac{1}{s+\frac{j-m}{2}}
+\sum_{j=0}^{\infty}C_j\frac{1}{s+\frac{j}{2}}+\sum_{j=0}^{\infty}
G_j\frac{1}{\left(s+\frac{j}{2}\right)^2}, \ s\to 0. 
\end{align}
In particular, $\zeta (s, \Delta_k^{\textup{rel/abs}})$ need not be regular at $s=0$. 
\end{cor}
The residue of the zeta-function at zero was computed by Cheeger in \cite{Che:SGOlv:AAS} and a crucial 
observation for the definition of analytic torsion in even dimensions has been made by Dar 
in \cite{Dar:IRT}
\begin{align}
\sum_{k=0}^{n+1}(-1)^k\, k \, \underset{s=0}{\textup{Res}}\, \zeta (s,  \Delta_k^{\textup{rel/abs}})=0.
\end{align}
We can now define the \emph{scalar analytic torsion} of $(\mathscr{C}(N),g)$ 
by differentiating at $s=0$ the following alternating weighted sum, rather than its individual summands
\begin{align}
\log T_{\textup{rel/abs}}(\mathscr{C}(N),E, g):=\left. \frac{d}{ds} \right|_{s=0}
\left(\frac{1}{2} \sum_{k=0}^{n+1}(-1)^k\, k \, \zeta_k(s,\Delta_k^{\textup{rel/abs}})\right).
\end{align}

As explicated in \cite{MV}, the scalar analytic torsions $T_{\textup{\textup{abs}}}(\mathscr{C}(N),g)$ and 
$T_{\textup{rel}}(\mathscr{C}(N),g)$ are related in even dimensions by Poincare duality 
\begin{align}\label{relabs}
\log T_{\textup{\textup{rel}}}(\mathscr{C}(N),E, g)=-\log T_{\textup{abs}}(\mathscr{C}(N),E, g).
\end{align}

The scalar analytic torsion is turned into the analytic torsion \emph{norm} by taking into account 
the de Rham cohomology. The determinant line on $(\mathscr{C}(N),E)$ is defined in terms of 
harmonic forms as
\begin{equation}\label{det-line}
\begin{split}
\det \mathscr{H}^*_{\textup{rel/abs}}(\mathscr{C}(N),E) &:=
\bigotimes_{k=0}^{\dim M}\left[\det \mathscr{H}^k_{\textup{rel/abs}}(\mathscr{C}(N),E)\right]^{(-1)^{k+1}}, \\
\det \mathscr{H}^k_{\textup{rel/abs}}(\mathscr{C}(N),E) &:=\bigwedge\nolimits^{\textup{top}}
\mathscr{H}^k_{\textup{rel/abs}}(\mathscr{C}(N),E). 
\end{split}
\end{equation}

\begin{defn}\label{t-norm-def}
The analytic torsion norm $\|\cdot \|^{RS, \, \textup{rel/abs}}_{(\mathscr{C}(N),E, g)}$ associated to $(\mathscr{C}(N),g)$ 
and the fixed $(E,\nabla, h^E)$, defined with respect to either relative or absolute boundary conditions, is the norm on 
$\det \mathscr{H}^*_{\textup{rel/abs}}(\mathscr{C}(N),E)$ given by
$$\|\cdot \|^{RS,  \, \textup{rel/abs}}_{(\mathscr{C}(N),E, g)}:=
T_{ \textup{rel/abs}}(\mathscr{C}(N), E, g)\|\cdot \|_{\det \mathscr{H}^*_{\textup{rel/abs}}(\mathscr{C}(N),E)},$$
where $\|\cdot \|_{\det \mathscr{H}^*_{\textup{rel/abs}}(\mathscr{C}(N),E)}$ is the norm 
on $\det \mathscr{H}^*_{\textup{rel/abs}}(\mathscr{C}(N),E)$ induced by the $L^2(g,h^E)$ 
norm on harmonic forms.
\end{defn}

\subsection{Statement of the Main Results}\label{main-results}

In this article we identify analytic torsion of $(\mathscr{C}(N),g)$, 
defined with respect to relative boundary conditions, in terms of cohomology and 
analytic torsion of the cone base $(N,g^N)$. We skip the index "rel" here, 
as we do not compare to analytic torsion with respect to absolute boundary conditions.

In the previous publication \cite{Ver:ATO}, we have 
derived an expression for the logarithm of analytic torsion on $(\mathscr{C}(N),g)$, 
comprised out of three types of contributions. 
\begin{thm}\label{BV-Theorem}\textup{(\!\!\cite{Ver:ATO}, Theorem 8.2)}
Let $\mathscr{C}(N)=(0,1]\times N, g=dx^2\oplus x^2g^N$ be an even-dimensional bounded cone over 
a closed oriented Riemannian manifold $(N,g^N)$ of dimension $n=\dim N$. Let $(E,\nabla, h^E)$ be a flat complex 
Hermitian vector bundle over $\mathscr{C}(N)$ and $(E_N,\nabla_N,h_N)$ its restriction to $N$.
Denote by $\Delta_{k,\ccl,N}$ the Laplacian on coclosed differential $k$-forms $\Omega^k_{\ccl}(N,E_N)$
and define 
\begin{align*}
&\A_k=\frac{(n-1)}{2}-k, \  F_k:=\{\nu = \sqrt{\eta +\A_k^2} \mid \eta \in \textup{Spec}\Delta_{k,\ccl,N}\backslash\{0\} \}, \\
&\zeta_{k,N}(s):=\sum_{\nu \in F_k} \nu^{-s}, \quad Re(s)\gg 0, \quad
 \delta_k:=\left\{ \begin{array}{cl} 1/2 & \textup{if} \ k=(n-1)/2, \\ 1 & \textup{otherwise}.\end{array}\right.
\end{align*}
Then the logarithm of the scalar analytic torsion of $(\mathscr{C}(N),g)$, defined with respect 
to relative boundary conditions, is given by a sum of a topological, spectral and the residual terms
$$\log T(\mathscr{C}(N),E,g)= \textup{Top}(N,E_N,g^N) + \textup{Tors}(N,E_N,g^N) + \textup{Res}(N,E_N,g^N),$$
where the topological term is an algebraic combination of Betti numbers 
$$\textup{Top}(N,E_N,g^N)=\sum_{k=0}^{(n-1)/2} \frac{(-1)^k}{2} \dim H^k(N,E_N)\log (n-2k+1).$$
The spectral term $\textup{Tors}(N,E_N,g^N)$ is expressed in terms of the scalar analytic torsion $T(N,E_N,g^N)$ of $(N,g^N)$ by
$$\textup{Tors}(N,E_N,g^N)=-\frac{1}{2}\log T(N,E_N,g^N).$$
The residual term is an intricate combination of residues of $\zeta_{k,N}(s)$
\begin{align*}
\textup{Res}(N,E_N,g^N) = \sum_{k=0}^{(n-1)/2}\frac{(-1)^k}{4} \delta_k\sum_{r=1}^{(n-1)/2} 
\underset{s=2r+1}{\textup{Res}}\, \zeta_{k,N}(s)\sum_{b=0}^{2r+1}A_{r,b}(\A_k)\frac{\Gamma'(b+r+1/2)}{\Gamma (b+r+1/2)},
\end{align*} 
where the coefficients $A_{r,b}(\A_k)$ are determined by certain recursive formulas, associated
to combinations of special functions.
\end{thm}

The computations in \cite{Ver:ATO} are performed in the untwisted setup, using the double summation 
method of Spreafico \cite{S, Spr:ZIF}.
However, any vector bundle $(E,\nabla,h^E)$ is product over $(\mathscr{C}(N)$, so the computations 
in \cite{Ver:ATO} carry over to the twisted setup along the same lines. Moreover, the formula 
for $\textup{Tors}(N,g^N)$ is presented in different terms in (\!\!\cite{Ver:ATO}, Theorem 8.2). 
Using (\!\!\cite{HS}, Lemma 7.2), an elegant combinatorial cancellation reduces the formula 
in (\!\!\cite{Ver:ATO}, Theorem 8.2) to the expression above. 

In view of our general results in \cite{Ver:ATO}, Melo-Hartmann-Spreafico evaluated later in \cite{HS, MHS} 
analytic torsion in the special case of a cone over $S^n$, identifying in even dimensions 
the singular contribution as the metric anomaly of the cone at its regular boundary by direct comparison 
of the expressions for analytic torsion and the metric anomaly in \cite{BruMa:AAF}. Due to the structure of the 
formulas, this approach via explicit comparison is limited in general case to lower dimensions. 

In this article we study this phenomena in the general case of an even dimensional bounded cone 
over any closed odd dimensional base manifold. Our main result here is that the third type of the contributions in 
the general formula of (\!\!\cite{Ver:ATO}, Theorem 8.2) indeed comes from the non-product metric structure 
of the cone at its regular boundary and vanishes if the cone metric is smoothened to a product 
away from the singularity. 

\begin{thm}\label{main1}
Let $(\mathscr{C}(N)=(0,1]\times N, g=dx^2\oplus x^2g^N$ be an even-dimensional bounded cone over 
a closed oriented Riemannian manifold $(N^n,g^N)$. Let $(E,\nabla, h^E)$ be a flat complex 
Hermitian vector bundle over $\mathscr{C}(N)$ and $(E_N,\nabla_N,h_N)$ its restriction to $N$.
Then the logarithm of the scalar analytic torsion of $(\mathscr{C}(N),g)$, defined with respect to the 
relative boundary conditions, is given by
\begin{align*}
\log T(\mathscr{C}(N),E,g)=& \sum_{k=0}^{(n-1)/2} \frac{(-1)^{k}}{2} \dim H^k(N,E_N)\log (n-2k+1) \\ 
- &\, \frac{1}{2}\log T(N,E_N,g^N) + \frac{\textup{rank(E)}}{2}\int_N B(\nabla^{TM}),
\end{align*}
where $T(N,E_N,g^N)$ is the Ray-Singer analytic torsion of $(N,g^N)$ and $B(\nabla^{TM})$ is a secondary characterstic 
class expressing the metric anomaly of $(\mathscr{C}(N),g)$ at $x=1$.
\end{thm}

\begin{cor}\label{main2}
Let $(\mathscr{C}(N)=(0,1]\times N, g=dx^2\oplus x^2g^N)$ be an even-dimensional bounded cone over 
a closed oriented Riemannian manifold $(N^n,g^N)$. Let $(E,\nabla, h^E)$ be a flat complex 
Hermitian vector bundle over $M$ and $(E_N,\nabla_N,h_N)$ its restriction to $N$. Let the metric $g_0$ on $M$ 
coincide with $g$ near the singularity at $x=0$ and be product $dx^2\oplus g^N$ in an open neighborhood of the 
boundary $\{1\}\times N$.  Then the Ray-Singer analytic torsion norm of $(\mathscr{C}(N),g_0)$ is given by
\begin{align*}
\log \|\cdot \|^{RS}_{(\mathscr{C}(N),E,g_0)} &= \sum_{k=0}^{(n-1)/2} \frac{(-1)^{k}}{2} \dim H^k(N,E_N) \log (n-2k+1) 
\\ &- \frac{1}{2}\log T(N,E_N,g^N) + \log \|\cdot \|_{\det \mathscr{H}^*(\mathscr{C}(N),E)},
\end{align*}
where $ \|\cdot \|_{\det \mathscr{H}^*(\mathscr{C}(N),E)}$ denotes the norm on the determinant line 
$\det \mathscr{H}^*(\mathscr{C}(N),E)$
induced by the $L^2(g,h^E)$-norm on square-integrable harmonic forms with relative boundary conditions.
\end{cor}

The main idea is to consider a cone-like cylinder, arising by truncating off the singularity from a cone.
Analytic torsion of such a cone-like cylinder $(\mathscr{C}_\e(N)=[\e,1]\times N, g=dx^2\oplus x^2g^N)$ 
with relative boundary conditions at $x=1$ and absolute boundary conditions at $x=\e$ is given by the 
analytic torsion of the exact cylinder plus twice the residual term $\textup{Res}(N,E_N,g^N)$. 
Application of \cite{BruMa:AAF} identifies the residual term as the metric anomaly at the regular boundary 
of the cone, leading to the results above. We have announced this approach in \cite{V2}, which 
has been followed later by Hartmann and Spreafico in \cite{HSpr:ZIF}.

The article is organized as follows. In Section \ref{section-torsion} we review the definition of analytic 
torsion on spaces with isolated conical singularities. In Section \ref{section-scaling} we discuss the 
metric anomaly formula for analytic torsion on even-dimensional manifolds with boundary, 
established by Br\"uning-Ma in \cite{BruMa:AAF}. In Section \ref{section-truncated} we compute analytic torsion 
of a cone-like cylinder, arising by truncating off the singularity from a cone. The main difference 
to the joint work with M\"uller \cite{MV} is the use of mixed boundary conditions rather than purely relative or absolute.
In Section \ref{section-anomaly} we apply this computation to identify the metric anomaly at the regular boundary 
of the cone and prove the main result. 

\section{Scaling Invariance of the Metric Anomaly for Analytic Torsion}
\label{section-scaling}

This section is the even-dimensional analogue of the corresponding discussion 
in \cite{MV}.  Denote by $\mathscr{A}\widehat{\otimes}\mathscr{B}$ 
the graded tensor product of $\Z_2-$graded algebras 
$\mathscr{A}, \mathscr{B}$. Put $\widehat{\mathscr{A}}:=I \widehat{\otimes} \mathscr{A}$ 
and identify $\mathscr{B}\equiv \mathscr{B} \widehat{\otimes} I$. 
Let $(X,g^X)$ be an even-dimensional compact oriented Riemannian manifold with boundary 
$\partial X$. Assume that over a collar neighborhood $\mathscr{U}\cong [0,1)\times \partial X,$
of the boundary $\partial X$, the Riemannian metric restricts to 
\begin{align*}
g^X\restriction \mathscr{U} = f(x) \left(dx^2 + g^{\partial X}\right), \ f\in C^{\infty}[0,1),
\end{align*}
This setup has been considered in (\!\!\cite{BruMa:AAF}, Examples 4.5). 
Let $R^{\partial X}$ be the curvature tensor of $(\partial X, g^{\partial X})$, defined in 
terms of the Levi-Civita connection of $g^{\partial X}$. Let $\{e_k\}_{k=1}^n$ 
be a local orthonormal frame on $(T\partial X,g^{\partial X})$ and $\{e^*_k\}_{k=1}^n$ the associated dual 
orthonormal frame on $T^*\partial X$. We denote by $\widehat{e^*_k}$ its canonical identification 
with elements in $\widehat{\Lambda T^*\partial X}$ and define
\begin{equation}\label{RS3}
\begin{split}
\dot{R}^{T\partial X}&= \frac{1}{2}\sum_{k,j}\langle e_k, R^{T\partial X} 
e_j\rangle e^*_k \wedge \widehat{e^*_j} \in 
\Lambda T^*\partial X \, \widehat{\otimes} \, \widehat{\Lambda T^*\partial X}\\
\dot{S}&=\frac{1}{4}f'(0) \sum_{k}e^*_k \wedge \widehat{e^*_k} 
\in \Lambda T^*\partial X \, \widehat{\otimes} \, \widehat{\Lambda T^*\partial X},
\end{split}
\end{equation}
For a general metric $g^X$, the two elements $\dot{R}^{T \partial X}$ and $\dot{S}$
are constructed in detail in (\!\!\cite{BruMa:AAF}, (1.15)), and the construction is more
intricate than the formulas above. However as above, $\dot{R}^{T\partial X}$ 
is defined in terms of the curvature of the Levi-Civita connection $\nabla^{TX}$; 
$\dot{S}$ measures the deviation from a metric product 
structure near the boundary. Moreover $\dot{R}^{T\partial X}$ and $\dot{S}^2$ are both homogeneous 
of degree two.

The linear map 
\begin{align}
\int^{B_{\partial X}}:\Lambda T^*\partial X \, \widehat{\otimes} \, 
\widehat{\Lambda T^*\partial X} \to \Lambda T^*\partial X,
\end{align}
is the Berezin integral, see (\!\!\cite{BruMa:AAF}, Section 1.1), which is non trivial only on elements, 
homogeneous of degree $\dim \partial X$. Then the secondary class $B(\nabla^{TX})$ 
is defined by the following expression, see (\!\!\cite{BruMa:AAF}, (1.17)) 
\begin{align}\label{RSpr:ZIF}
B(\nabla^{TX})=- \int_0^1 \frac{du}{u}\int^{B_{\partial X}} \exp \left(-\frac{1}{2}\dot{R}^{T\partial X}
-u^2\dot{S}^2\right) \sum_{k=1}^{\infty} \frac{(u\dot{S})^k}{2\Gamma(k/2+1)},
\end{align}
Let $(E,\nabla^E)$ be flat complex vector bundle over $X$, equipped with
the associated flat Hermitian metric $h^E$. Assume, $\partial X=Y_1\cup Y_2$ is a 
disjoint union of two closed manifolds and we define the analytic torsion norm 
with respect to absolute boundary conditions at $Y_1$ and relative boundary conditions at $Y_2$.
If $g^X_i$, $i=1,2$, are two Riemannian metrics on $X$, then (\!\!\cite{BruMa:AAF2}, Theorem 3.4) 
asserts for the corresponding Ray-Singer 
metrics $\|\cdot \|^{RS}_{(X,E;g^X_i)}$, $i=1,2$ on $\det \mathscr{H}^*(X,E)$ 
\begin{align}\label{BM-thm}
\log \left(\frac{\|\cdot \|^{RS}_{(X,E;g^X_1)}}{\|\cdot \|^{RS}_{(X,E;g^X_2)}}
\right)=\frac{\textup{rank}(E)}{2}
\left(\int_{Y_1}-\int_{Y_2}\right)\left[ B(\nabla^{TX}_2)-B(\nabla^{TX}_1)\right].
\end{align}

The sign difference to (\!\!\cite{BruMa:AAF2}, Theorem 3.4)  is due to different determinant line conventions. 
An important property of the secondary class $B(\nabla^{TX})$ is its scaling invariance.
\begin{prop}\label{Scaling}
Let $g^X_1$ and $g^X_2=s\cdot g^X_1, s>0,$ be any pair of Riemannian metrics on a compact manifold 
$(X,\partial X)$ with boundary and let $\nabla^{TX}_1, \nabla^{TX}_2$ denote the corresponding 
Levi-Civita connections. Then 
\begin{align}
 B(\nabla^{TX}_1)=B(\nabla^{TX}_2).
\end{align}
\end{prop}

\begin{proof}
Denote the constituents of \eqref{RSpr:ZIF}, defined by the Riemannian metrics $g^X_1$ and $g^X_2$ by 
$(\dot{R}^{T\partial X}_1,\dot{S}_1)$ and $(\dot{R}^{T\partial X}_2,\dot{S}_2)$ respectively. 
By definition in (\!\!\cite{BruMa:AAF}, (1.15)) we find
\begin{equation} 
\dot{R}^{\partial T}(g^X_2)=s\dot{R}^{\partial T}(g^X_1), 
\quad \dot{S}(g^X_2)=\sqrt{s}\dot{S}(g^X_1).
\end{equation}
The Berezin integral also depends on the metric, which we denote by $\int^{B_{\partial X}}_j,j=1,2$. 
By definition, see (\!\!\cite{BruMa:AAF}, (1.1)), we find
\begin{align}
\int^{B_{\partial X}}_2=s^{-\dim \partial X /2} \int^{B_{\partial X}}_1.
\end{align}
Since $\dot{R}^{T\partial X}$ and $\dot{S}^2$ are both homogeneous of degree two, and the Berezin integral 
is non-trivial only on terms homogeneous of degree $\dim \partial X$, we conclude

\begin{align*}
B(\nabla^{TX}_2)&=- \int_0^1 \frac{du}{u}\int^{B_{\partial X}}_2 \exp \left(-\frac{1}{2}\dot{R}^{T\partial X}_2
-u^2\dot{S}^2_2\right) \sum_{k=1}^{\infty} \frac{(u\dot{S}_2)^k}{2\Gamma(k/2+1)}\\
&= - \int_0^1 \frac{du}{u}\int^{B_{\partial X}}_1 \exp \left(-\frac{1}{2}\dot{R}^{T\partial X}_2
-u^2\dot{S}^2_2\right) \sum_{k=1}^{\infty} \frac{(u\dot{S}_2)^k}{2\Gamma(k/2+1)} \cdot s^{-\dim \partial X /2}\\
&= - \int_0^1 \frac{du}{u}\int^{B_{\partial X}}_1 \exp \left(-\frac{1}{2}\dot{R}^{T\partial X}_1
-u^2\dot{S}^2_1\right) \sum_{k=1}^{\infty} \frac{(u\dot{S}_1)^k}{2\Gamma(k/2+1)}=B(\nabla^{TX}_1).
\end{align*} 
\end{proof}

\section{Analytic Torsion of the Cone-like Cylinder}
\label{section-truncated}

\subsection{Decomposition of the de Rham Complex of the Cone-like Cylinder}\label{decomposition}

This section is the even-dimensional analogue of the corresponding discussion 
in \cite{MV}. We repeat the argument here once again for completeness. 
Consider the product manifold $\mathscr{C}_I (N)=I\times N$
over an even-dimensional closed Riemannian manifold $(N^n,g^N), n=\dim N$, 
where either $I=(0,1)$ modelling the cone, or $I=[\epsilon, 1], \epsilon >0$ modelling a cone-like cylinder. 
The Riemannian metric on $\mathscr{C}_I(N)$ is given by a warped product 
$$g = dx^2 \oplus x^2g^N, x\in I.$$ 
Let $(E,\nabla, h^E)$ be a flat complex Hermitian vector bundle over $\mathscr{C}_I(N)$ and $(E_N,\nabla_N,h_N)$
its restriction to the cross-section $N$. 
Consider the associated twisted de Rham complex $(\Omega^*_0(\mathscr{C}_I(N),E),\nabla_*)$, 
where $\Omega^*_0(\mathscr{C}_I(N),E)$ are the $E$-valued differential forms with compact support in 
$\mathscr{C}_I(N)$, and $\nabla_*$ denotes the covariant differential induced by the flat connection $\nabla$ 
on $E$. Consider the separation of variables map
\begin{align}\label{separation}
\Psi_k : C^{\infty}_0(I,\Omega^{k-1}(N,E_N)\oplus \Omega^k(N,E_N))\to \Omega_0^k(\mathscr{C}_I(N),E) \\
(\w_{k-1},\w_k)\mapsto x^{k-1-n/2}\w_{k-1}\wedge dx + x^{k-n/2}\w_k, \nonumber
\end{align}
where $\w_k,\w_{k-1}$ are identified with their pullback to $\mathscr{C}_I(N)$ under the 
projection $\pi: I\times N\to N$ onto the second factor, and $x$ is the on $I\subset \R^+$. 
The map $\Psi_k$ extends to an isometry, cf. \cite{BruSee:AIT}
\begin{align*}
\Psi_k: L^2(I, L^2(\Omega^{k-1}(N,E_N)\oplus \Omega^k(N,E_N), g^N, h_N), dx)\to 
L^2(\Omega_0^k(\mathscr{C}_I(N),E), g,h^E).
\end{align*}
We continue under the unitary transformation above henceforth and obtain under that identification, as in (\!\!\cite{BruSee:AIT}, (5.5))
\begin{equation}\label{derivative} 
\begin{split}
\nabla_k&\equiv \Psi_{k+1}^{-1} \nabla_k \Psi_k= \left( \begin{array}{cc}0&(-1)^k\partial_x\\0&0\end{array}\right)+\frac{1}{x}
\left( \begin{array}{cc}\nabla_{k-1,N}&c_k\\0&\nabla_{k,N}\end{array}\right), \\
\nabla_k^*&\equiv\Psi_k^{-1} \nabla_k^* \Psi_{k+1}= \left( \begin{array}{cc}0&0\\(-1)^{k+1}\partial_x&0\end{array}\right)+
\frac{1}{x}\left( \begin{array}{cc}\nabla_{k-1,N}^*&0\\c_k&\nabla_{k,N}^*\end{array}\right), 
\end{split}
\end{equation}
where $\nabla_{k,N}$ is the de Rham differential on $\Omega^k(N,E_N)$, and 
$$c_k:=(-1)^k\left(k-\frac{n}{2}\right).$$
Following \cite{L:P}, 
we decompose the de Rham complex $(\Omega^*_0(\mathscr{C}_I(N),E),\nabla_*)$ into a direct sum of subcomplexes of two types. 
The first type of the subcomplexes is given as follows. Let $\psi \in \Omega^k(N,E_N)$ be a coclosed 
eigenform of the Laplacian $\Delta_{k,N}$ on $\Omega^k(N,E_N)$ with eigenvalue $\eta >0$. 
We consider the following four associated pairs
\begin{equation}\label{xi1}
\begin{split}
&\xi_1:=(0,\psi)\in \Omega^{k-1}(N,E_N)\oplus \Omega^{k}(N,E_N), \\ &\xi_2:=(\psi,0)\in 
\Omega^{k}(N,E_N)\oplus \Omega^{k+1}(N,E_N), \\
&\xi_3:=(0,\nabla_N\psi/\sqrt{\eta})\in \Omega^{k}(N,E_N)\oplus \Omega^{k+1}(N,E_N), \\ 
&\xi_4:=(\nabla_N\psi/\sqrt{\eta},0)\in \Omega^{k+1}(N,E_N)\oplus \Omega^{k+2}(N,E_N).
\end{split}
\end{equation}
Denote by $\langle \xi_1,\xi_2,\xi_3,\xi_4\rangle$ the vector space, spanned by the four vectors in \eqref{xi1}.
$C^{\infty}_0(I,\langle \xi_1,\xi_2,\xi_3,\xi_4\rangle)$ is invariant under $\nabla,\nabla^*$ and we obtain a subcomplex
\begin{align}\label{complex-1}
0 \rightarrow C_0^{\infty}(I,\left< \xi_1\right>) \xrightarrow{\nabla^{\psi}_0} C_0^{\infty}
(I,\left<\xi_2,\xi_3\right>) \xrightarrow{\nabla^{\psi}_1}C_0^{\infty}(I,\left<\xi_4\right>)
\rightarrow 0,
\end{align}
where $\nabla^{\psi}_0,\nabla^{\psi}_1$ take the following form with respect to the chosen basis:
\begin{align*}
\nabla_0^{\psi}=\left(\begin{array}{c}(-1)^k\partial_x+\frac{c_k}{x}\\ x^{-1}\sqrt{\eta}\end{array}\right), 
\quad \nabla_1^{\psi}=\left(x^{-1}\sqrt{\eta}, \ (-1)^{k+1}\partial_x+\frac{c_{k+1}}{x}\right).
\end{align*}
The associated Laplacians are of the following form
\begin{align}\label{psi-laplacians}
\Delta_0^{\psi}:=(\nabla_0^{\psi})^*\nabla_0^{\psi}=-\partial_x^2+\frac{1}{x^2}\left[\eta +
\left(k+\frac{1}{2}-\frac{n}{2}\right)^2-\frac{1}{4}\right]=\nabla_1^{\psi}(\nabla_1^{\psi})^*=:
\Delta_2^{\psi}.
\end{align}
under the identification of any $\phi =f\cdot \xi_{i}\in C^{\infty}_0(I,\langle \xi_{i}\rangle),i=1,4$ 
with its scalar part $f \in C^{\infty}_0(I)$. We continue under this identification from here on. 
Subcomplexes \eqref{complex-1} always come in pairs on oriented cones. The twin subcomplex is 
constructed by considering $\phi:=*_N\psi \in \Omega^{n-k}(N,E_N)$. Then $\nabla^*_N\phi/\sqrt{\eta}$ is 
again a coclosed eigenform of the Laplacian on $\Omega^{n-k-1}(N,E_N)$ with eigenvalue $\eta$, and we put 
\begin{equation}\label{xi2}
\begin{split}
&\widetilde{\xi_1}:=(0,\nabla^*_N\phi/\sqrt{\eta})\in \Omega^{n-k-2}(N,E_N)\oplus \Omega^{n-k-1}(N,E_N), 
\\ &\widetilde{\xi_2}:=(\nabla^*_N\phi/\sqrt{\eta},0)\in \Omega^{n-k-1}(N,E_N)\oplus \Omega^{n-k}(N,E_N), \\
&\widetilde{\xi_3}:=(0,\phi)\in \Omega^{n-k-1}(N,E_N)\oplus \Omega^{n-k}(N,E_N), \\ &\widetilde{\xi_4}
:=(\phi,0)\in \Omega^{n-k}(N,E_N)\oplus \Omega^{n-k+1}(N,E_N).
\end{split}
\end{equation}
Denote by $\langle \widetilde{\xi_1},\widetilde{\xi_2},\widetilde{\xi_3},\widetilde{\xi_4}\rangle$ 
the vector space, spanned by the four vectors in \eqref{xi2}.
$C_0^{\infty}(I,\langle \widetilde{\xi_1},\widetilde{\xi_2},\widetilde{\xi_3},\widetilde{\xi_4}\rangle)$ 
is invariant under the action of $\nabla,\nabla^*$ and we obtain a subcomplex
\begin{align}\label{complex-2}
0 \rightarrow C_0^{\infty}(I,\langle \widetilde{\xi_1}\rangle) \xrightarrow{\nabla_0^{\phi}} 
C_0^{\infty}(I,\langle\widetilde{\xi_2},\widetilde{\xi_3}\rangle) \xrightarrow{\nabla_1^{\phi}}
C_0^{\infty}(I,\langle\widetilde{\xi_4}\rangle)\rightarrow 0.
\end{align}
Computing explicitly the action of the exterior derivative \eqref{derivative} on the basis elements $\widetilde{\xi_i}$ we find
\begin{align*}
\nabla_0^{\phi}=\left(\begin{array}{c}(-1)^{n-k-1}\partial_x+\frac{c_{n-k-1}}{x}\\ x^{-1}\sqrt{\eta}\end{array}\right), 
\quad \nabla_1^{\phi}=\left(x^{-1}\sqrt{\eta}, \ (-1)^{n-k}\partial_x+\frac{c_{n-k}}{x}\right).
\end{align*}
As before we compute the corresponding Laplacians and find
\begin{align}\label{phi-laplacians}
\Delta_0^{\phi}=\Delta_2^{\phi}=-\partial_x^2+\frac{1}{x^2}\left[\eta +\left(k+\frac{1}{2}-
\frac{n}{2}\right)^2-\frac{1}{4}\right]=\Delta_0^{\psi}=\Delta_2^{\psi},
\end{align}
where the operators are again identified with their scalar actions.

The second type of the subcomplexes comes from the harmonic differential forms 
$H^k(N,E_N)$ over the base manifold $N$. Fix an orthonormal basis $\{u_i\}$ of $H^k(N,E_N)$
and observe that any subspace $C^{\infty}_0(I,\langle0\oplus u_i,u_i\oplus 0\rangle)$ 
is invariant under $\nabla,\nabla^*$. Consequently we obtain a subcomplex of the de Rham complex
\begin{equation}\label{complex-3}\begin{split}
0\to C^{\infty}_0(I,\langle 0\, \oplus \, &u^k_i\rangle)\xrightarrow{\nabla^H_k} 
C^{\infty}_0(I,\langle u^k_i\oplus 0\rangle) \to 0, \\
&\nabla^H_k=(-1)^k\partial_x+\frac{c_k}{x},
\end{split}
\end{equation}
where the action of $\nabla^H_k$ is identified with its scalar action, as before. The Laplacians of the complex are given by
\begin{equation}
\begin{split}
H^k_0:&=(\nabla^H_k)^*\nabla^H_k=-\partial_x^2+\frac{1}{x^2}
\left(\left(\frac{(n-1)}{2}-k\right)^2-\frac{1}{4}\right), \\ 
H^k_1:&=\nabla^H_k(\nabla^H_k)^*=-\partial_x^2+\frac{1}{x^2}
\left(\left(\frac{(n+1)}{2}-k\right)^2-\frac{1}{4}\right).
\end{split}
\end{equation}  

\subsection{The Relative and Absolute Boundary Conditions}

The covariant differential $\nabla_k$ is a priori defined on the domain $\Omega^k_0(\mathscr{C}_I(N),E)$ 
of differential forms with compact support. 
Let $\nabla_{k,\min}$ denote the graph closure of $\nabla_k$  and $\nabla_{k,\max}$ its maximal closed extension in $L^2(\Omega_0^*(\mathscr{C}_I(N),E),g,h^E)$.  
Despite the de Rham differential $\nabla_{k}$ not being elliptic, there is still 
a well-defined trace on $\dom (\nabla_{k,\max})$ by 
the trace theorem of Paquet in \cite{P}.
\begin{thm}\label{trace-theorem} \textup{(\!\!\cite{P}, Theorem 1.9)} 
Let $(X,g^X)$ be a compact oriented Riemannian manifold, possibly with isolated conical singularities, 
and with smooth boundary $\partial X$. Let $\iota: \partial X \hookrightarrow X$ be the natural inclusion. 
Let $(E,\nabla, h^E)$ be a flat complex Hermitian vector bundle over $X$ and $(E_{\partial X},\nabla_{\partial X},h_{\partial X})$
its restriction to the boundary, again a flat complex Hermitian vector bundle over $\partial X$.
Then the pullback $\iota^*:\Omega^k(X,E) \to \Omega^k(\partial X, E_{\partial X})$ with 
$\Omega^k(\partial X, E_{\partial X})=\{0\}$ for $k=\dim X$, extends continuously to a linear surjective map 
\begin{align}
\iota^*:\dom (\nabla_{k,\max})\rightarrow H^{-1/2}(\nabla_{k,\partial X}),
\end{align} 
where $\nabla_{k,\partial X}$ is the de Rham differential on $\Omega^k(\partial X, E_{\partial X})$,  
$H^{-1/2}(\Omega^k(\partial X, E_{\partial X}))$ the $(-1/2)$-th Sobolev space on $\partial X$ and
\begin{align*}
 H^{-1/2}(\nabla_{k,\partial X}):=\{\w \in H^{-1/2}(\Omega^k(\partial X, E_{\partial X})) \mid 
\nabla_{k,\partial X} \w \in H^{-1/2}(\Omega^{k+1}(\partial X, E_{\partial X}))\},
\end{align*}
is a Hilbert space under the obvious graph-norm.
\end{thm}

\begin{remark}
The trace theorem (\!\!\cite{P}, Theorem 1.9) is stated in the untwisted setup on compact (non-singular) 
Riemannian manifolds. Extension to flat Hermitian vector bundles is straightforward. Moreover, the analysis 
localizes to an open neighborhood of the boundary $\partial X$, so the trace theorem carries over to compact 
Riemannian manifolds with singular structure 
away from $\partial X$.
\end{remark}

Recall the relative and the absolute self-adjoint extensions of the Laplacian
\begin{equation}
\begin{split}
&\Delta_k^{\textup{rel}}:=\nabla^*_{k,\min}\nabla_{k,\min}+\nabla_{k-1,\min}
\nabla^*_{k-1,\min}, \\
&\Delta_k^{\textup{abs}}:=\nabla^*_{k,\max}\nabla_{k,\max}+\nabla_{k-1,\max}
\nabla^*_{k-1,\max}. 
\end{split}
\end{equation}
By definition of the maximal and minimal closed extension, $\nabla^*_{k,\max}=\nabla^t_{k,\min}$ and 
consequently $\dom (\nabla_{m-k-1,\min})=*\dom (\nabla^*_{k,\max})$.
Hence, Theorem \ref{trace-theorem} implies 
\begin{equation}
\begin{split}
&\dom (\nabla_{k,\min})\subseteq \{\w\in \dom (\nabla_{k,\max})|\iota^*\w=0\}, \\
&\dom (\nabla^*_{k,\max})\subseteq \{\w\in \dom (\nabla^t_{k,\max})|\iota^*(*\w)=0\},\\
&\dom (\Delta^{\textup{rel}}_k)\subseteq \{\w \in \dom (\Delta_{k,\max})| \iota^*\w=0, \iota^*(\nabla^t_{k-1}\w)=0\}, \\
&\dom (\Delta^{\textup{abs}}_k)\subseteq \{\w \in \dom (\Delta_{k,\max})| \iota^*(*\w)=0, \iota^*(*\nabla_{k}\w)=0\}.
\end{split}
\end{equation}
By the Hodge decomposition of $\Omega^*(N,E_N)$, the de Rham complex 
$(\Omega_0^*(\mathscr{C}_I(N), E), \nabla_*)$ 
decomposes completely into subcomplexes of the three types \eqref{complex-1}, \eqref{complex-2} and \eqref{complex-3}. 
It has been observed in (\!\!\cite{Ver:ATO}, Theorem 3.5) that in each degree $k$ this induces a compatible decomposition 
for the relative and absolute extension of the Laplacian. In the classical language of 
\cite{W2} we have a decomposition into reducing subspaces of the Laplacians. 
Hence the Laplacians $\Delta_k^{\textup{rel}}, \Delta_k^{\textup{abs}}$ induce self-adjoint extensions of the 
Laplacians $\Delta_{j}^{\psi}, \Delta_{j}^{\phi},j=0,2,$ and $H^k_0,H^k_1$. 

In case of $I=[\e,1], \e>0,$ it becomes necessary to deal with mixed boundary conditions for 
the Laplace operator on $\mathscr{C}_I(N)$. Let $\gamma \in C^{\infty}_0(\e,1]$ be a 
cut-off function, vanishing identically near $\{x=\e\}$ and being identically one near $\{x=1\}$. 
Define a mixed closed extension $D_k$ of the exterior derivative $\nabla_k$, and the associated mixed self adjoint 
extension of $\Delta_k$ by
\begin{equation}\label{rel-bc}
\begin{split}
\dom (D_k)&:=\{\w \in \dom (\nabla_{k,\max}) \mid \gamma \w \in  \dom (\nabla_{k,\min})\}, \\
\Delta_{k}^{\textup{mix}}&:=D_k^*D_k + D_{k-1}D_{k-1}^*.
\end{split}
\end{equation}
By construction, $\Delta_k^{\textup{mix}}$ has relative boundary conditions at $\{x=1\}$ and absolute boundary conditions 
at $\{x=\e\}$. In case of $I=(0,1]$ the usual relative self adjoint extension of the Laplace operator shall still be 
denoted by $\Delta_k^{\textup{mix}}$. The decomposition of the de Rham complex into 
subcomplexes of the three types \eqref{complex-1}, \eqref{complex-2} and \eqref{complex-3} 
induces a compatible decomposition for the mixed self-adjoint extension of the Laplacian, by similar arguments 
as in (\!\!\cite{Ver:ATO}, Theorem 3.5). For linguistic convenience we refer to the relative self-adjoint extension on the 
cone as the mixed extension again, so we don't need to switch between two notions hence and forth. 

In order to discuss the mixed boundary conditions explicitly, note that by the classical 
theory of linear differential equations for any 
element $f$ of $\dom (\Delta_{j,\max}^{\psi}),\dom (\Delta_{j,\max}^{\phi}),j=0,2,$ or $\mathscr{D}(H^k_{i,\max}),i=0,1,$ 
$f$ and its derivative $f'$ are both locally absolutely continuous in $I$, with well-defined values at $x\in \partial I$, more 
precisely at $x=1$ in case $I=(0,1]$, and $x\in \{\e,1\}$ in case $I=[\e,1]$. 
Hence the following boundary conditions are well-defined 
\begin{equation*}
B_N^k(x)f:=f'(x)+(-1)^{k+1}c_{k}\frac{f(x)}{x}, \quad B_D(x)f:=f(x), \quad x\in \partial I.
\end{equation*}
In case of $I=(0,1]$ boundary conditions at $x=0$ need to be posed. 
By the well-known analysis, compare \cite{BruSee:AIT}, \cite{Che:SGOlv:AAS}, 
see also an overview \cite{Ver:ZDF}, 
any solution $f\in L^2(0,1)$ to 
\begin{align}
 -\frac{d^2f}{dx^2}+\frac{1}{x^2}\left(\nu^2-\frac{1}{4}\right)f=g\in L^2(0,1),
\end{align}
admits an asymptotic expansion of the form
\begin{equation}
 f(x)\sim \left\{ 
\begin{array}{ll}
c_1(f)\sqrt{x}+c_2(f)\sqrt{x}\log(x)+O(x^{3/2}),  &\nu=0, \\
c_1(f)x^{\nu+1/2}+c_2(f)x^{-\nu+1/2}+O(x^{3/2}),  &\nu\in (0,1),\\
O(x^{3/2}), &\nu\geq 1,
\end{array}
\right. \ x\to 0,
\end{equation}
where the coefficients $c_1(f)$ and $c_2(f)$ depend only on $f$. Consequently the following boundary conditions at $x=0$ 
are well-defined
\begin{align}
 B_N(0)f:=
\left\{ 
\begin{array}{ll} c_1(f), & \nu \in [0,1), \\ 0, &\nu\geq 1,\end{array}\right.
\quad 
B_D(0)f:=
\left\{ 
\begin{array}{ll} c_2(f), & \nu \in [0,1), \\ 0, &\nu\geq 1.\end{array}\right.
\end{align}

\begin{prop}\label{rel-bc-prop}
Let $(\mathscr{C}_I(N)=I\times N, g = dx^2 \oplus x^2g^N)$ be a cone-type manifold over a closed even-dimensional Riemannian 
manifold $(N^n,g^N)$ with either $I=(0,1]$ or $I=[\e, 1], \e>0$. Consider the Laplacians $\Delta_{j}^{\psi},\Delta_{j}^{\phi}, j=0,2,$ 
of the subcomplex-pair \eqref{complex-1} and \eqref{complex-2}, and the Laplacians $H^k_0,H^k_1,$ of the subcomplex \eqref{complex-3}. 
The domains of their mixed self-adjoint extensions 
are given as follows. For $I=[\e,1]$
\begin{align*}
&\mathscr{D}(\Delta_{0,\textup{mix}}^{\psi})=\{f\in \mathscr{D}(\Delta_{0,\max}^{\psi})\mid B_N^{n-k}(\e)f=0, \ B_D(1)f=0\}, \\
&\mathscr{D}(\Delta_{0,\textup{mix}}^{\phi})=\{f\in \mathscr{D}(\Delta_{0,\max}^{\phi})\mid B_N^{k+1}(\e)f=0, \ B_D(1)f=0\}, \\
&\mathscr{D}(\Delta_{2,\textup{mix}}^{\psi})=\{f\in \mathscr{D}(\Delta_{2,\max}^{\psi})\mid B_D(\e)f=0, \ B_N^{k+1}(1)f=0\}, \\
&\mathscr{D}(\Delta_{2,\textup{mix}}^{\phi})=\{f\in \mathscr{D}(\Delta_{2,\max}^{\phi})\mid B_D(\e)f=0, \ B_N^{n-k}(1)f=0\}, \\
&\mathscr{D}(H^k_{1,\textup{mix}})=\{f\in \mathscr{D}(H^k_{1,\max})\mid B_D(\e)f=0, \ B_N^k(1)f=0\}, \\
&\mathscr{D}(H^k_{0,\textup{mix}})=\{f\in \mathscr{D}(H^k_{0,\max})\mid B_N^{n-k}(\e)f=0, \ B_D(1)f=0\}.
\end{align*}
For $I=(0,1]$ the domains are given by
\begin{align*}
&\mathscr{D}(\Delta_{0,\textup{mix}}^{\psi,\, \phi})=\{f\in\mathscr{D}(\Delta_{0,\max}^{\psi,\, \phi}) \mid B_D(0)f=0, \ B_D(1)f=0\}, \\ 
&\mathscr{D}(\Delta_{2,\textup{mix}}^{\psi})=\{f\in \mathscr{D}(\Delta_{2,\max}^{\psi})\mid B_D(0)f=0, \ B_N^{k+1}(1)f=0\}, \\
&\mathscr{D}(\Delta_{2,\textup{mix}}^{\phi})=\{f\in \mathscr{D}(\Delta_{2,\max}^{\phi})\mid B_D(0)f=0, \ B_N^{n-k}(1)f=0\}, \\
&\mathscr{D}(H^k_{1,\textup{mix}})=\{f\in \mathscr{D}(H^k_{1,\max})\mid B_N(0)f=0, \ B_N^k(1)f=0\}, \\
&\mathscr{D}(H^k_{0,\textup{mix}})=\{f\in \mathscr{D}(H^k_{0,\max})\mid B_D(0)f=0, \ B_D(1)f=0\}.
\end{align*}
\end{prop}

\begin{proof}
Boundary conditions at $x=1$ in case $I=(0,1]$, at $x\in \{\e,1\}$ in case $I=[\e,1]$, 
follow for the individual mixed self-adjoint extensions from \eqref{rel-bc}, the explicit form of the de Rham differentials 
\eqref{derivative} and the fact that for any $x\in \partial I$ 
and the inclusion $\iota_x:\{x\}\times N \hookrightarrow \mathscr{C}_I(N)$, we have 
$\iota^*_x(f_{k-1},f_k)=f_k(x)$ for any $(f_{k-1},f_k)\in \dom (\Delta_{\max})$ 
with $f_k$ continuous at $x$. Boundary conditions at $x=0$ in case $I=(0,1]$ have been determined in 
(\!\!\cite{Ver:ZDF}, Corollary 2.14) and (\!\!\cite{Ver:ATO}, Proposition 3.6 and 3.7).
\end{proof}

\subsection{Difference of Analytic Torsions of Cone-like Cylinder and the Cone}\label{difference}

We now distinguish between the cases $I=(0,1]$ and $I=[\e,1]$; 
and write $(\mathscr{C}(N)=(0,1]\times N, g=dx^2\oplus x^2g^N)$ 
for the bounded cone over an odd-dimensional closed Riemannian manifold $(N^n,g^N)$, 
and $(\mathscr{C}_\e (N)=[\e,1]\times N, g=g\restriction \mathscr{C}_\e (N))$ 
for its truncation. Let us write $\Delta_k^{\textup{mix}}$ and $\Delta_{k,\e}^{\textup{mix}}$ for the Laplacians 
with mixed boundary conditions on $k$-forms associated to $(\mathscr{C}(N),g)$ and 
$(\mathscr{C}_\e (N),g)$, respectively. Put
\begin{align}
 T(\e,s):=\frac{1}{2} \, \sum_{k=0}^{\dim \mathscr{C}(N)} (-1)^k\cdot k 
\cdot \left(\zeta(s,\Delta_{k,\e}^{\textup{\textup{mix}}})-\zeta(s,\Delta_{k}^{\textup{\textup{mix}}})\right).
\end{align}
$T(\e,s)$ is related to the scalar analytic torsions of $(\mathscr{C}(N),g)$ and 
$(\mathscr{C}_\e (N),g)$ by
\begin{align}
T'(\e,0)=\log T_{\textup{mix}}(\mathscr{C}_\e (N), E, g)-
\log T_{\textup{rel}}(\mathscr{C}(N), E, g).
\end{align}
Consider the decomposition of the de Rham complex in Section \ref{decomposition}. 
For each fixed degree $k$, the subcomplexes \eqref{complex-1} and \eqref{complex-2} are determined 
by a coclosed eigenform $\psi \in \Omega^k(N,E_N)$ of the Laplacian $\Delta_{k,N}$ with eigenvalue $\eta>0$. 
Denote the $\eta$-dependence by writing $\psi\equiv \psi(\eta)$, with eigenvalues $\eta$ coming from the set
$$E_k:=\textup{Spec}(\Delta_{k,\ccl,N})\backslash \{0\}.$$
Mixed boundary conditions for the Laplacians $\Delta_{j}^{\psi(\eta)},\Delta_{j}^{\phi(\eta)}, j=0,2,$ 
of the subcomplex-pair \eqref{complex-1} and \eqref{complex-2}, and 
the Laplacians $H^k_0,H^k_1,$ of the subcomplex \eqref{complex-3} 
are discussed in Proposition \ref{rel-bc-prop}. Here we distinguish operators on 
$(\mathscr{C}_\e (N),g)$ by an additional $\e$-subscript. 

\begin{defn} \label{zetas} Put for $\textup{Re}(s)\gg 0$
\begin{equation}\begin{split}
\zeta_{k,H}(s,\e)&:= \dim H^k(N,E_N) 
\left(\zeta(s, H^k_{0,\e,\textup{\textup{mix}}}) - \zeta(s, H^k_{0,\textup{\textup{mix}}})\right), \\
\zeta_k(s,\e)&:=\sum_{\eta \in E_k}\left(\zeta(s, \Delta^{\psi(\eta)}_{2,\e,\textup{mix}})+ 
\zeta(s, \Delta^{\phi(\eta)}_{2,\e,\textup{mix}}) -\zeta(s, \Delta^{\psi(\eta)}_{0,\e, \textup{mix}}) - 
\zeta(s, \Delta^{\phi(\eta)}_{0,\e,\textup{mix}}) \right)\\ 
&- \sum_{\eta \in E_k}\left(\zeta(s, \Delta^{\psi(\eta)}_{2,\textup{mix}})+ 
\zeta(s, \Delta^{\phi(\eta)}_{2,\textup{mix}}) -\zeta(s, \Delta^{\psi(\eta)}_{0, \textup{mix}}) - 
\zeta(s, \Delta^{\phi(\eta)}_{0,\textup{mix}}) \right).
\end{split}
\end{equation}
\end{defn}

$\zeta_{k,H}(s,\e)$ and $\zeta_k(s,\e)$ contribute to $T(\e,s)$, cf. (\!\!\cite{Ver:ATO}, (4.3), (4.4)), as follows
\begin{align}\label{T}
T(\e,s)=\frac{1}{2}\sum_{k=0}^{(n-1)/2-1}(-1)^k\, \delta_k \, \zeta_k(s,\e) + 
\frac{1}{2}\sum_{k=0}^{n}(-1)^{k+1}\zeta_{k,H}(s,\e), \ \textup{Re}(s)\gg 0,
\end{align}
where
$$ 
\delta_k:=\left\{ \begin{array}{cl} 1/2 & \textup{if} \ k=(n-1)/2, \\ 1 & \textup{otherwise}.
\end{array}\right.
$$
Evaluation of $\zeta_k'(0,\e)$ requires application of the double summation method, 
introduced by Spreafico in \cite{S}, \cite{Spr:ZIF} and applied by the us to 
derive the general formula for analytic torsion of a bounded cone in 
\cite{Ver:ATO}, see Theorem \ref{BV-Theorem}.
Evaluation of $\zeta'_{k,H}(0,\e)$ reduces to an explicit computation of finitely 
many zeta-determinants and application of \cite{Les:DOR}.
We begin with the evaluation of $\zeta_k'(0,\e)$ for each fixed degree $k$ along 
the lines of (\!\!\cite{Ver:ATO}, Section 6).

\begin{prop}\label{N-prop}
Let the contour 
$\Lambda_c:=\{\lambda \in \C| |\textup{arg}(\lambda -c)|=\pi /4\}$
be oriented counter-clockwise, for any $c>0$.
Fix the branch of logarithm in $\C\backslash \R^-$. Put 
\begin{align*}
\A_k:=\frac{(n-1)}{2}-k, \quad 
\nu(\eta):=\sqrt{\eta + \A_k^2}, \ \eta \in \textup{Spec}\Delta_{k,\ccl,N}\backslash \{0\}
\end{align*}
Let $c(\eta)=c_0/(2\nu(\eta)^2),$ where $c_0>0$ is a fixed positive number, 
smaller than the lowest non-zero eigenvalue of $\Delta^{\textup{\textup{mix}}}_*$ and 
$\Delta^{\textup{\textup{mix}}}_{*,\epsilon}$, such that $c(\eta)<1$ for all $\eta\in E_k$.
Then $\zeta_k(s,\e)$ 
admits the following integral representation for 
$\textup{Re}(s)\gg 0$
\begin{align*}
\zeta_k(s,\e)=\sum_{\eta \in E_k} \nu(\eta)^{-2s} 
\frac{s^2}{\Gamma(s+1)}\int_0^{\infty}\frac{t^{s-1}}{2\pi i}
\int_{\wedge_{c(\eta)}}\frac{e^{-\lambda t}}{-\lambda}\, 
t_{\eta, \e}^{k}(\lambda)  d\lambda \, dt,
\end{align*}
with $t_{\eta,\e}^{k}(\lambda)$ given in terms of zeta-determinants
\begin{equation}\label{tm}
\begin{split}
t_{\nu,\e}^{k}(\lambda)= -& \log \frac{\det_{\zeta} \left(\Delta^{\psi(\eta)}_{2,\e,\textup{mix}}+\nu^2z^2\right)}{\det \left(\Delta^{\psi(\eta)}_{2,\e,\textup{mix}}\right)}
-\log \frac{\det_{\zeta} \left(\Delta^{\phi(\eta)}_{2,\e,\textup{mix}}+\nu^2z^2\right)}{\det \left(\Delta^{\phi(\eta)}_{2,\e,\textup{mix}}\right)} \\
+& \log \frac{\det_{\zeta} \left(\Delta^{\psi(\eta)}_{0,\e,\textup{mix}}+\nu^2z^2\right)}{\det \left(\Delta^{\psi(\eta)}_{0,\e,\textup{mix}}\right)}
+\log \frac{\det_{\zeta} \left(\Delta^{\phi(\eta)}_{0,\e,\textup{mix}}+\nu^2z^2\right)}{\det \left(\Delta^{\phi(\eta)}_{0,\e,\textup{mix}}\right)} 
\\
 +& \log \frac{\det_{\zeta} \left(\Delta^{\psi(\eta)}_{2,\textup{mix}}+\nu^2z^2\right)}{\det \left(\Delta^{\psi(\eta)}_{2,\textup{mix}}\right)}
+\log \frac{\det_{\zeta} \left(\Delta^{\phi(\eta)}_{2,\textup{mix}}+\nu^2z^2\right)}{\det \left(\Delta^{\phi(\eta)}_{2,\textup{mix}}\right)} \\
-& \log \frac{\det_{\zeta} \left(\Delta^{\psi(\eta)}_{0,\textup{mix}}+\nu^2z^2\right)}{\det \left(\Delta^{\psi(\eta)}_{0,\textup{mix}}\right)}
-\log \frac{\det_{\zeta} \left(\Delta^{\phi(\eta)}_{0,\textup{mix}}+\nu^2z^2\right)}{\det \left(\Delta^{\phi(\eta)}_{0,\textup{mix}}\right)} .
\end{split}
\end{equation}
\end{prop}
\begin{proof}
Recall that the spectrum entering $\zeta_k(s,\epsilon)$ is the union of spectra for the Laplacians 
$\Delta^{\psi(\eta)}_{2,\e,\textup{\textup{mix}}}, \Delta^{\psi(\eta)}_{2,\textup{\textup{mix}}}$ and 
$\Delta^{\phi(\eta)}_{2,\e,\textup{\textup{mix}}}, \Delta^{\phi(\eta)}_{2,\textup{\textup{mix}}}$, where $\eta$ runs over $E_k$. 
For any choice 
$$
L(\eta)\in \left\{\Delta^{\psi(\eta)}_{2,\e,\textup{\textup{mix}}}, \Delta^{\psi(\eta)}_{2,\textup{\textup{mix}}}, 
\Delta^{\phi(\eta)}_{2,\e,\textup{\textup{mix}}}, \Delta^{\phi(\eta)}_{2,\textup{\textup{mix}}}\right\}, \eta\in E_k,
$$
the spectrum $\textup{Spec}\, L(\eta)\subset \R^+$ is strictly positive. Indeed, $\textup{Spec}\, L(\eta)$ is contained in
the spectrum of the non-negative Laplace operator on the truncated or full cone, 
and its zero eigenvalues arise in both cases only from harmonic forms $H^*(N,E_N)$. 
Resolvent of $L(\eta)$ is trace class, cf. \cite{Les:DOR}, and from 
Definition \ref{zetas} we infer for $\textup{Re}(s)\gg 0$
\begin{align}\label{integral}
\zeta_k(s,\e)=\sum_{\eta\in E_k} \nu(\eta)^{-2s} \frac{1}{\Gamma(s)}\int_0^{\infty}t^{s-1}\frac{1}{2\pi i}
\int_{\wedge_{c(\eta)}}e^{-\lambda t}h_{\eta,\e}^{k}(\lambda) \,   d\lambda dt,
\end{align}
where
\begin{align*}
h_{\eta,\e}^{k}(\lambda) &=
   \textup{Tr}\left(\lambda - \nu(\eta)^{-2}\Delta^{\psi(\eta)}_{2,\epsilon, \textup{\textup{mix}}}\right)^{-1} 
+ \textup{Tr}\left(\lambda - \nu(\eta)^{-2}\Delta^{\phi(\eta)}_{2,\epsilon, \textup{\textup{mix}}}\right)^{-1} \\
&- \textup{Tr}\left(\lambda - \nu(\eta)^{-2}\Delta^{\psi(\eta)}_{0,\epsilon, \textup{\textup{mix}}}\right)^{-1} 
- \textup{Tr}\left(\lambda - \nu(\eta)^{-2}\Delta^{\phi(\eta)}_{0,\epsilon, \textup{\textup{mix}}}\right)^{-1} \\
&-  \textup{Tr}\left(\lambda - \nu(\eta)^{-2}\Delta^{\psi(\eta)}_{2,\textup{\textup{mix}}}\right)^{-1} 
- \textup{Tr}\left(\lambda - \nu(\eta)^{-2}\Delta^{\phi(\eta)}_{2, \textup{\textup{mix}}}\right)^{-1}\\
&+ \textup{Tr}\left(\lambda - \nu(\eta)^{-2}\Delta^{\psi(\eta)}_{0,\textup{\textup{mix}}}\right)^{-1} 
+ \textup{Tr}\left(\lambda - \nu(\eta)^{-2}\Delta^{\phi(\eta)}_{0, \textup{\textup{mix}}}\right)^{-1}.
\end{align*}
For any choice of 
$$
L(\eta)\in \left\{\Delta^{\psi(\eta)}_{2,\e,\textup{\textup{mix}}}, \Delta^{\psi(\eta)}_{2,\textup{\textup{mix}}}, 
\Delta^{\phi(\eta)}_{2,\e,\textup{\textup{mix}}}, \Delta^{\phi(\eta)}_{2,\textup{\textup{mix}}}\right\}, \eta\in E_k,
$$
we find by (\!\!\cite{Les:DOR}, Proposition 4.6) that, enumerating $\textup{Spec}\, L(\eta)=\{\lambda_i\}_{i=1}^{\infty}$
in increasing order, the series 
\begin{equation}\label{det-det}
\log \frac{\det_{\zeta} (L(\eta)-\nu(\eta)^{2}\lambda)}{\det_{\zeta} L(\eta)} = \sum_{i=1}^{\infty} \log \left(1-\frac{\nu(\eta)^{2}\lambda}{\lambda_i}\right).
\end{equation}
converges and by the choice of the logarithm branch is holomorphic in 
$\lambda \in \C \backslash \{x\in \R \mid x > c(\eta)\}$. Moreover, 
\begin{align}\label{tr-det}
\textup{Tr} \left(\frac{L(\eta)}{\nu(\eta)^2}-\lambda\right)^{-1} =-\frac{d}{d\lambda} \log \frac{\det_{\zeta} 
(L(\eta)-\nu(\eta)^{2}\lambda)}{\det_{\zeta} L(\eta)}.
\end{align} 
By the definition of $c(\eta)>0$, \eqref{det-det} is holomorphic in an open neighborhood of 
the contour $\Lambda_{c(\eta)}$, and so we may integrate \eqref{integral} by parts 
first in $\lambda$ then in $t$, and obtain
\begin{align}
\zeta_k(s,\e)&= \sum_{\eta\in E_k} \nu(\eta)^{-2s}
\frac{1}{\Gamma(s)}\int_0^{\infty}t^{s-1}\frac{1}{2\pi i}
 \int_{\wedge_{c(\eta)}}e^{-\lambda t}h_{\eta,\e}^{k}(\lambda)  d\lambda dt 
\\&= \sum_{\eta\in E_k} \nu(\eta)^{-2s} 
\frac{s^2}{\Gamma(s+1)}\int_0^{\infty}t^{s-1}\frac{1}{2\pi i}
\int_{\wedge_{c(\eta)}}\frac{e^{-\lambda t}}{-\lambda}t_{\eta,\e}^{k}(\lambda) d\lambda dt,
\end{align}
where 
\begin{equation}
\begin{split}
t_{\nu,\e}^{k}(\lambda)= -& \log \frac{\det_{\zeta} \left(\Delta^{\psi(\eta)}_{2,\e,\textup{mix}}+\nu^2z^2\right)}
{\det \left(\Delta^{\psi(\eta)}_{2,\e,\textup{mix}}\right)}
-\log \frac{\det_{\zeta} \left(\Delta^{\phi(\eta)}_{2,\e,\textup{mix}}+\nu^2z^2\right)}{\det \left(\Delta^{\phi(\eta)}_{2,\e,\textup{mix}}\right)} \\
+& \log \frac{\det_{\zeta} \left(\Delta^{\psi(\eta)}_{0,\e,\textup{mix}}+\nu^2z^2\right)}{\det \left(\Delta^{\psi(\eta)}_{0,\e,\textup{mix}}\right)}
+\log \frac{\det_{\zeta} \left(\Delta^{\phi(\eta)}_{0,\e,\textup{mix}}+\nu^2z^2\right)}{\det \left(\Delta^{\phi(\eta)}_{0,\e,\textup{mix}}\right)} 
\\
 +& \log \frac{\det_{\zeta} \left(\Delta^{\psi(\eta)}_{2,\textup{mix}}+\nu^2z^2\right)}{\det \left(\Delta^{\psi(\eta)}_{2,\textup{mix}}\right)}
+\log \frac{\det_{\zeta} \left(\Delta^{\phi(\eta)}_{2,\textup{mix}}+\nu^2z^2\right)}{\det \left(\Delta^{\phi(\eta)}_{2,\textup{mix}}\right)} \\
-& \log \frac{\det_{\zeta} \left(\Delta^{\psi(\eta)}_{0,\textup{mix}}+\nu^2z^2\right)}{\det \left(\Delta^{\psi(\eta)}_{0,\textup{mix}}\right)}
-\log \frac{\det_{\zeta} \left(\Delta^{\phi(\eta)}_{0,\textup{mix}}+\nu^2z^2\right)}{\det \left(\Delta^{\phi(\eta)}_{0,\textup{mix}}\right)} .
\end{split}
\end{equation}
\end{proof}  

\begin{lemma}\label{det-Bessel}
For any $\nu>0$ and $z \in \C$
\begin{equation*}
\begin{split}
&\frac{\det_{\zeta} \left(\Delta^{\psi(\eta)}_{2, \textup{\textup{mix}}}+\nu^2z^2\right)}
{\det \left(\Delta^{\psi(\eta)}_{2,\textup{\textup{mix}}}\right)}= \frac{2^{\nu}\Gamma(\nu)}
{(\nu z)^{\nu}(1+\A_k/\nu)}\left(\nu z I'_{\nu}(\nu z)+\A_k I_{\nu}(\nu z)\right),\\ 
&\frac{\det_{\zeta} \left(\Delta^{\phi(\eta)}_{2,\textup{\textup{mix}}}+\nu^2z^2\right)}
{\det \left(\Delta^{\phi(\eta)}_{2,\textup{\textup{mix}}}\right)}=\frac{2^{\nu}\Gamma(\nu)}
{(\nu z)^{\nu}(1-\A_k/\nu)}\left(\nu z I'_{\nu}(\nu z)-\A_k I_{\nu}(\nu z)\right), \\
&\frac{\det_{\zeta} \left(\Delta^{\psi(\eta)}_{0, \textup{mix}}+\nu^2z^2\right)}
{\det \left(\Delta^{\psi(\eta)}_{0,\textup{mix}}\right)}= \frac{\det_{\zeta} 
\left(\Delta^{\phi(\eta)}_{0, \textup{mix}}+\nu^2z^2\right)}{\det \left(\Delta^{\phi(\eta)}_{0,\textup{mix}}\right)}=
\frac{2^{\nu}\Gamma(\nu+1)}{(\nu z)^{\nu}}I_{\nu}(\nu z). 
\end{split}
\end{equation*}

\begin{equation*}
\begin{split}
\frac{\det_{\zeta} \left(\Delta^{\psi(\eta)}_{2,\e,\textup{\textup{mix}}}+\nu^2z^2\right)}
{\det \left(\Delta^{\psi(\eta)}_{2,\e,\textup{\textup{mix}}}\right)}&=  
2\nu \, \frac{(\nu z I'_\nu(\nu z) + \A_k I_\nu (\nu z) )K_\nu(\nu z \e)}{(\nu + \A_k)\e^{-\nu}+ (\nu -\A_k)\e^\nu} \\
& \times
\left(1-\frac{ \nu z K'_{\nu}(\nu z)+\A_k K_{\nu}(\nu z) }
{ \nu z I'_{\nu}(\nu z)+\A_k I_{\nu}(\nu z)}\cdot \frac{I_{\nu}(\nu z \e)}{K_{\nu}(\nu z \e)}\right),
\\
\frac{\det_{\zeta} \left(\Delta^{\phi(\eta)}_{2,\e,\textup{\textup{mix}}}+\nu^2z^2\right)}
{\det \left(\Delta^{\phi(\eta)}_{2,\e,\textup{\textup{mix}}}\right)}&=
2\nu \, \frac{(\nu z I'_\nu(\nu z) - \A_k I_\nu (\nu z) )K_\nu(\nu z \e)}{(\nu - \A_k)\e^{-\nu}+ (\nu + \A_k)\e^\nu} \\
& \times
\left(1-\frac{ \nu z K'_{\nu}(\nu z)-\A_k K_{\nu}(\nu z) }
{ \nu z I'_{\nu}(\nu z)-\A_k I_{\nu}(\nu z)}\cdot \frac{I_{\nu}(\nu z \e)}{K_{\nu}(\nu z \e)}\right),
\\
\frac{\det_{\zeta} \left(\Delta^{\psi(\eta)}_{0,\e,\textup{\textup{mix}}}+\nu^2z^2\right)}
{\det \left(\Delta^{\psi(\eta)}_{2,\e,\textup{\textup{mix}}}\right)}&= 
2\nu \, \frac{(- \nu z \e K'_\nu(\nu z \e) + \A_k K_\nu (\nu z \e) ) I_\nu(\nu z)}{(\nu + \A_k)\e^{-\nu}+ (\nu - \A_k)\e^\nu} \\
& \times
\left(1-\frac{K_{\nu}(\nu z)}{I_{\nu}(\nu z)}\cdot \frac{ \nu z \e I'_{\nu}(\nu z\e)
- \A_k I_{\nu}(\nu z \e)}{ \nu z\e K'_{\nu}(\nu z\e) - \A_k K_{\nu}(\nu z\e)}\right),
\\
\frac{\det_{\zeta} \left(\Delta^{\phi(\eta)}_{0,\e,\textup{\textup{mix}}}+\nu^2z^2\right)}
{\det \left(\Delta^{\phi(\eta)}_{2,\e,\textup{\textup{mix}}}\right)}&=
2\nu \, \frac{(- \nu z \e K'_\nu(\nu z \e) - \A_k K_\nu (\nu z \e) ) I_\nu(\nu z)}{(\nu - \A_k)\e^{-\nu}+ (\nu + \A_k)\e^\nu} \\
& \times
\left(1-\frac{K_{\nu}(\nu z)}{I_{\nu}(\nu z)}\cdot \frac{ \nu z \e I'_{\nu}(\nu z\e)
+ \A_k I_{\nu}(\nu z \e)}{ \nu z\e K'_{\nu}(\nu z\e) + \A_k K_{\nu}(\nu z\e)}\right).
\end{split}
\end{equation*}
\end{lemma}

\begin{proof}
We evaluate zeta-determinants using their explicit relation with the 
normalized solutions of the operators, established by Lesch in (\!\!\cite{Les:DOR}, Theorem 1.2). 
The first four equations have been evaluated in (\!\!\cite{Ver:ATO}, Corollary 6.3)
\begin{equation}\label{det1}
\begin{split}
&\frac{\det_{\zeta} \left(\Delta^{\psi(\eta)}_{2, \textup{\textup{mix}}}+\nu^2z^2\right)}
{\det \left(\Delta^{\psi(\eta)}_{2,\textup{\textup{mix}}}\right)}= \frac{2^{\nu}\Gamma(\nu)}
{(\nu z)^{\nu}(1+\A_k/\nu)}\left(\nu z I'_{\nu}(\nu z)+\A_k I_{\nu}(\nu z)\right),\\ 
&\frac{\det_{\zeta} \left(\Delta^{\phi(\eta)}_{2,\textup{\textup{mix}}}+\nu^2z^2\right)}
{\det \left(\Delta^{\phi(\eta)}_{2,\textup{\textup{mix}}}\right)}=\frac{2^{\nu}\Gamma(\nu)}
{(\nu z)^{\nu}(1-\A_k/\nu)}\left(\nu z I'_{\nu}(\nu z)-\A_k I_{\nu}(\nu z)\right), \\
&\frac{\det_{\zeta} \left(\Delta^{\psi(\eta)}_{0, \textup{mix}}+\nu^2z^2\right)}
{\det \left(\Delta^{\psi(\eta)}_{0,\textup{mix}}\right)}= \frac{\det_{\zeta} 
\left(\Delta^{\phi(\eta)}_{0, \textup{mix}}+\nu^2z^2\right)}{\det \left(\Delta^{\phi(\eta)}_{0,\textup{mix}}\right)}=
\frac{2^{\nu}\Gamma(\nu+1)}{(\nu z)^{\nu}}I_{\nu}(\nu z).
\end{split}
\end{equation}
In order to evaluate zeta determinants of 
$\Delta^{\psi(\eta)}_{2,\e,\textup{\textup{mix}}}, \Delta^{\phi(\eta)}_{2,\e,\textup{\textup{mix}}}$ and 
$\Delta^{\psi(\eta)}_{0,\e,\textup{\textup{mix}}}, \Delta^{\phi(\eta)}_{0,\e,\textup{\textup{mix}}}$, 
consider solutions $f_{\psi, \nu}(\cdot, z)$ and $f_{\phi,\nu}(\cdot, z)$
of $(\Delta^{\psi(\eta)}_{2,\e,\textup{\textup{mix}}}+z^2)$ and $(\Delta^{\phi(\eta)}_{2,\e,\textup{\textup{mix}}}+z^2)$, 
respectively, normalized at $x=1$. By definition, see (\!\!\cite{Les:DOR}, (1.38a), (1.38b)), these are solutions of the respective operators, 
satisfying relative boundary conditions at $x=1$ and normalized by  $f_{\psi, \nu}(1, z)=1$ and $f_{\phi,\nu}(1, z)=1$, i.e. 
\begin{align*}
\begin{array}{lll}
(\Delta^{\psi(\eta)}_{2,\e}+z^2)f_{\psi, \nu}(\cdot, z)=0, & 
f'_{\psi, \nu}(1, z) + (-1)^kc_{k+1}f_{\psi, \nu}(1, z)=0, & f_{\psi, \nu}(\cdot, z)=1, \\ 
(\Delta^{\phi(\eta)}_{2,\e}+z^2)f_{\phi, \nu}(1, z)=0, 
& f'_{\phi, \nu}(1, z) + (-1)^{n-k+1}c_{n-k}f_{\phi, \nu}(1, z)=0, & f_{\phi, \nu}(1, z)=1.
\end{array}
\end{align*}
Normalized solutions are uniquely determined and explicit computations lead to the following expressions
\begin{equation}\label{f-mu}
\begin{split}
f_{\psi, \nu}(x, z)&=(zI'_{\nu}(z)+\A_kI_{\nu}(z))\sqrt{x}K_{\nu}(zx) - (zK'_{\nu}(z)+\A_kK_{\nu}(z))\sqrt{x}I_{\nu}(zx), \\
f_{\phi, \nu}(x, z)&=(zI'_{\nu}(z)-\A_kI_{\nu}(z))\sqrt{x}K_{\nu}(zx) - (zK'_{\nu}(z)-\A_kK_{\nu}(z))\sqrt{x}I_{\nu}(zx), \\
f_{\psi, \nu}(x, 0)&=\frac{1}{2\nu} (\nu-\A_k)x^{\nu+1/2} + \frac{1}{2\nu} (\nu +\A_k)x^{-\nu+1/2}, \\
f_{\phi, \nu}(x, 0)&=\frac{1}{2\nu} (\nu+\A_k)x^{\nu+1/2} + \frac{1}{2\nu} (\nu -\A_k)x^{-\nu+1/2},
\end{split}
\end{equation}
where we use 
\begin{align}
K_{\nu}(z)I'_{\nu}(z)-K'_{\nu}(z)I_{\nu}(z)=\frac{1}{z}.
\end{align}
In view of (\!\!\cite{Les:DOR}, Theorem 1.2), this fact being due to Burghelea-Friedlander-Kappeler in \cite{BFK:OTD} 
in the non-singular setup,  we find
\begin{equation}\label{det2}
\begin{split}
\frac{\det_{\zeta} \left(\Delta^{\psi(\eta)}_{2,\e,\textup{\textup{mix}}}+\nu^2z^2\right)}
{\det \left(\Delta^{\psi(\eta)}_{2,\e,\textup{\textup{mix}}}\right)}=
\frac{f_{\psi, \nu}(\e, \nu z)}{f_{\psi, \nu}(\e, 0)}, \\
\frac{\det_{\zeta} \left(\Delta^{\phi(\eta)}_{2,\e,\textup{\textup{mix}}}+\nu^2z^2\right)}
{\det \left(\Delta^{\phi(\eta)}_{2,\e,\textup{\textup{mix}}}\right)}=
\frac{f_{\phi, \nu}(\e, \nu z)}{f_{\phi, \nu}(\e, 0)}.
\end{split}
\end{equation}
Plugging in the expressions \eqref{f-mu} we arrive at the statemement of the lemma. 
The determinants associated to $\Delta^{\psi(\eta)}_{0,\e,\textup{\textup{mix}}}, 
\Delta^{\phi(\eta)}_{0,\e,\textup{\textup{mix}}}$ are discussed along the same lines.
\end{proof}

In particular, applying Lemma \ref{det-Bessel} several cancellations lead to a 
representation of $t_{\eta,\e}^{k}(\lambda)$ in terms of Bessel functions with
$\nu\equiv \nu(\eta)$ and $z=\sqrt{-\lambda}$, where we use the fixed branch of logarithm in $\C\backslash \R^-$, 
extended by continuity to one of the sides of the cut
\begin{equation}\label{t-bessel}
\begin{split}
t_{\nu,\e}^{k}(\lambda)=&-2 \log K_{\nu}(\nu z \e) -\log \left(1-\frac{\A_k^2}{\nu^2}\right) - 2\log (\nu)\\
+&\log \left(-\nu z \e K'_{\nu}(\nu z \e) + \A_k K_{\nu}(\nu z \e) \right)\\
+&\log \left(-\nu z \e K'_{\nu}(\nu z \e) - \A_k K_{\nu}(\nu z \e)\right)\\
- &\log \left(1-\frac{ \nu z K'_{\nu}(\nu z)+\A_k K_{\nu}(\nu z)}
{ \nu z I'_{\nu}(\nu z)+\A_k I_{\nu}(\nu z)}\cdot \frac{I_{\nu}(\nu z \e)}{K_{\nu}(\nu z \e)}\right) \\ 
- &\log \left(1-\frac{ \nu z K'_{\nu}(\nu z)-\A_k K_{\nu}(\nu z)}{ \nu z I'_{\nu}(\nu z)-\A_k I_{\nu}(\nu z)}
\cdot \frac{I_{\nu}(\nu z \e)}{K_{\nu}(\nu z \e)} \right) \\
+ &\log \left(1-\frac{K_{\nu}(\nu z)}{I_{\nu}(\nu z)}\cdot \frac{ \nu z \e I'_{\nu}(\nu z\e)
+ \A_k I_{\nu}(\nu z \e)}{ \nu z\e K'_{\nu}(\nu z\e) + \A_k K_{\nu}(\nu z\e)}\right) \\
+ &\log \left(1-\frac{K_{\nu}(\nu z)}{I_{\nu}(\nu z)}\cdot \frac{ \nu z \e I'_{\nu}(\nu z\e)
- \A_k I_{\nu}(\nu z \e)}{ \nu z\e K'_{\nu}(\nu z\e) - \A_k K_{\nu}(\nu z\e)}\right)
\end{split}
\end{equation}

For the arguments below we need to summarize some facts on Bessel functions.
Consider expansions of Bessel-functions for large arguments and fixed order, cf. (\!\!\cite{AbrSte:HOM}, p.377).
 For the modified Bessel functions of first kind we have
\begin{equation}\label{large-arg-I}
\begin{split}
I_{\nu}(z)=\frac{e^z}{\sqrt{2\pi z}}\left(1+O\left(\frac{1}{z}\right)\right), \\
 I'_{\nu}(z)=\frac{e^z}{\sqrt{2\pi z}}\left(1+O\left(\frac{1}{z}\right)\right), 
\end{split}
\quad |z|\to \infty.
\end{equation}
Expansions for modified Bessel functions of second kind are
\begin{equation}\label{large-arg-K}
\begin{split}
K_{\nu}(z)=\sqrt{\frac{\pi}{2 z}}e^{-z}\left(1+O\left(\frac{1}{z}\right)\right), \\
K'_{\nu}(z)=-\sqrt{\frac{\pi}{2 z}}e^{-z}\left(1+O\left(\frac{1}{z}\right)\right),
\end{split}
\quad |z|\to \infty.
\end{equation}
$|\textup{arg}(z)|<\pi /2$ is region of validity for \eqref{large-arg-I} and \eqref{large-arg-I}, and the expansions
in particular hold for $z=\sqrt{-\lambda}$ with $\lambda \in \Lambda_c$ large. 
For small arguments and positive orders $\nu>0$ we have the following expansions
\begin{equation}\label{small}
\begin{split}
I_{\nu}(z)\sim \frac{z^{\nu}}{2^{\nu}\Gamma (\nu+1)}, \quad K_{\nu}(z)\sim 2^{\nu-1}\frac{\Gamma(\nu)}{z^{\nu}}&, \\ 
I'_{\nu}(z)\sim \frac{z^{\nu-1}}{2^{\nu}\Gamma(\nu)}, \quad K'_{\nu}(z)\sim-2^{\nu-1}\frac{\Gamma(\nu+1)}{z^{\nu+1}}&,
\end{split}
\quad \textup{as} \ |z|\to 0.
\end{equation}
Consider expansions of Bessel-functions for large order $\nu>0$, cf. (\!\!\cite{Olv:AAS}, Section 7). 
For any $\{z\in \C||\textup{arg}(z)|<\pi /2\}\cup \{ix|x\in (-1,1)\},$ put $t:=(1+z^2)^{-1/2}$ and 
$\xi:=1/t+\log (z/(1+1/t))$. For the modified Bessel functions of first kind we then have
\begin{equation}\label{large-nu-I}
\begin{split}
I_{\nu}(\nu z)  = \frac{1}{\sqrt{2\pi \nu}}\frac{e^{\nu \xi}}{(1+z^2)^{1/4}} 
\left[1+\sum_{r=1}^{N-1}\frac{u_r(t)}{\nu^r} + \frac{\eta_{N,1}(\nu, z)}{\nu^{N}} \right], & \\
I'_{\nu}(\nu z)  = \frac{1}{\sqrt{2\pi \nu}
} \frac{e^{\nu \xi}}{z (1+z^2)^{-1/4}}
\left[1+\sum_{r=1}^{N-1}\frac{v_r(t)}{\nu^r} + \frac{\eta_{N,2}(\nu, z)}{\nu^{N}} \right]. &
\end{split}
\end{equation} 
Expansions for modified Bessel functions of second kind are
\begin{equation}\label{large-nu-K}
\begin{split}
K_{\nu}(\nu z) = \sqrt{\frac{\pi}{2 \nu}} \frac{e^{-\nu \xi}}{(1+z^2)^{1/4}} 
\left[1+\sum_{r=1}^{N-1}\frac{u_r(t)}{(-\nu)^r} + \frac{\eta_{N,3}(\nu, z)}{(-\nu)^{N}} \right], & \\
K'_{\nu}(\nu z) =  -\sqrt{\frac{\pi}{2 \nu}} \frac{e^{-\nu \xi}}{z (1+z^2)^{-1/4}}
\left[1+\sum_{r=1}^{N-1}\frac{v_r(t)}{(-\nu)^r} + \frac{\eta_{N,4}(\nu, z)}{(-\nu)^{N}} \right]. &
\end{split}
\end{equation}
The error terms $\eta_{N,i}(\nu, z)$ are bounded for large $\nu$ 
uniformly in any compact subset of $\{z\in \C||\textup{arg}(z)|<\pi /2\}\cup \{ix|x\in (-1,1)\}$, 
see the analysis of validity regions for the expansions 
\eqref{large-nu-I} and \eqref{large-nu-K} in (\!\!\cite{Olv:AAS}, Section 8). 
For $\lambda \in \Lambda_c$ with any $0<c<1$, the induced $z=\sqrt{-\lambda}$ is contained in that region of validity, 
where we use the fixed branch of logarithm in $\C\backslash \R^-$, extended by continuity to one of the sides of the cut. The 
coefficients $u_r(t),v_r(t)$ are polynomial in $t$ and defined via a recursive relation, see (\!\!\cite{Olv:AAS}, (7.10)).

As in (\!\!\cite{BKD:HKA}, (3.15)) we have for any fixed $\A\in \R$ 
the following expansion as $\nu\to \infty$ in terms of orders
\begin{equation}\label{polynom2}
\begin{split}
&\log  \left(1+\sum_{r=1}^{N}\frac{u_r(t)}{(\pm\nu)^r}+ O(\nu^{-N-1})\right) 
\sim \sum_{r=1}^{\infty}\frac{D_r(t)}{(\pm\nu)^r} + O(\nu^{-N-1}), \\
&\log \left[ \left(1+\sum_{k=1}^{N}\frac{v_r(t)}{(\pm\nu)^r} \right)
+ \frac{\A}{(\pm\nu)}t\left(1+\sum_{r=1}^{N-1}\frac{u_r(t)}{(\pm\nu)^r}\right)+ O(\nu^{-N-1})\right] 
 \\ &\sim \sum_{r=1}^{N}\frac{M_r(t, \A)}{(\pm\nu)^r} + O(\nu^{-N-1}),
\end{split}
\end{equation}
where by the polynomial structure of $u_r(t)$ and $v_r(t)$, the coefficients 
$D_r(t)$ and $M_r(t,\A)$ are polynomial in $t$, see also (\!\!\cite{BKD:HKA}, (3.7), (3.16)), with 
\begin{align}\label{MD-polynom}
D_r(t)=\sum_{b=0}^{r}x_{r,b}t^{r+2b}, \quad 
M_r(t, \A)=\sum_{b=0}^{r}z_{r,b}(\A)t^{r+2b}.
\end{align}
As a consequence of (\!\!\cite{BGKE:ZFD}, (4.24))
\begin{align}\label{DM}
M_{r}(1,\A) = D_{r}(1)-\frac{(-\A)^{r}}{r}.
\end{align}

\begin{prop}\label{contour-change}
There exist $\e,c>0$ small enough such that for 
$\textup{Re}(s)\gg 0$
\begin{align*}
\zeta_k(s,\e)=
\frac{s^2}{\Gamma(s+1)}\int_0^{\infty}\frac{t^{s-1}}{2\pi i}
\int_{\wedge_{c}}\frac{e^{-\lambda t}}{-\lambda}\, 
T^k_{\e}(s,\lambda) d\lambda \, dt, \quad 
T^k_{\e}(s,\lambda) =\sum_{\eta \in E_k} t_{\eta, \e}^{k}(\lambda)  \, \nu(\eta)^{-2s} 
\end{align*}
\end{prop}

\begin{proof}
Consider the representation of $t_{\eta, \e}^{k}(\lambda)$ in \eqref{t-bessel}
in terms of Bessel functions. We need to investigate its behaviour for large $\eta$, 
or equivalently for large $\nu(\eta)$.
Let $\{z\in \C||\textup{arg}(z)|<\pi /2\}\cup \{ix|x\in (-1,1)\}$ 
and $t_{\e}:=(1+(\e z)^2)^{-1/2}$. By \eqref{large-nu-K} we find
\begin{equation}\label{log1}
\begin{split}
\log &\left(-\nu z \e K'_{\nu}(\nu z \e)+\A_k K_{\nu}(\nu z \e)\right) +
\log \left(-\nu z \e K'_{\nu}(\nu z \e)-\A_k K_{\nu}(\nu z \e)\right) \\
= \, &\log \left[ \left(1+\sum_{k=1}^{N-1}\frac{v_r(t_{\e})}{(-\nu)^r} \right)+ \frac{\A_k}
{(-\nu)}t_{\e}\left(1+\sum_{r=1}^{N-2}\frac{u_r(t_{\e})}{(-\nu)^r}\right) + 
\frac{\kappa_{N,1}(\nu,z\e)}{(-\nu)^N}\right] \\
&+\log \left[ \left(1+\sum_{k=1}^{N-1}\frac{v_r(t_{\e})}{(-\nu)^r} \right)- 
\frac{\A_k}{(-\nu)}t_{\e}\left(1+\sum_{r=1}^{N-2}\frac{u_r(t_{\e})}{(-\nu)^r}\right) +
\frac{\kappa_{N,2}(\nu,z\e)}{(-\nu)^N}\right] ,
\end{split}
\end{equation}
where the error terms
\begin{equation}
\begin{split}
\kappa_{N,1}(\nu,z\e)&=
\eta_{N,4}(\nu,z\e) + (\A_kt_{\e}) \eta_{N-1,3}(\nu,z\e) \\
\kappa_{N,2}(\nu,z\e)&=
\eta_{N,4}(\nu,z\e) - (\A_kt_{\e}) \eta_{N-1,3}(\nu,z\e)
\end{split}
\end{equation}
are bounded for large $\nu$ 
uniformly in any compact subset of $\{z\in \C||\textup{arg}(z)|<\pi /2\}\cup \{ix|x\in (-1,1)\}$.
Employing \eqref{large-nu-I} and \eqref{large-nu-K} we find with
$\xi:=1/t+\log (z/(1+1/t))$ and $\xi_{\e}:=1/t_{\e}+\log (\e z/(1+1/t_{\e}))$
\begin{equation}\label{log2}
\begin{split}
\log \left(1-\frac{ \nu z K'_{\nu}(\nu z)+\A_k K_{\nu}(\nu z)}
{ \nu z I'_{\nu}(\nu z)+\A_k I_{\nu}(\nu z)}\cdot \frac{I_{\nu}(\nu z \e)}{K_{\nu}(\nu z \e)}\right) 
=\log \left( 1- e^{2\nu(\xi_{\e}-\xi)}(1+\kappa_1(\nu,z))\right),
\end{split}
\end{equation}
where the error term $\kappa(\nu,z)$ is again bounded for large $\nu$ 
uniformly in any compact subset of $\{z\in \C||\textup{arg}(z)|<\pi /2\}\cup \{ix|x\in (-1,1)\}$. 
Similarly, 
\begin{equation}\label{log3}
\begin{split}
\log \left(1-\frac{ \nu z K'_{\nu}(\nu z)-\A_k K_{\nu}(\nu z)}{ \nu z I'_{\nu}(\nu z)-\A_k I_{\nu}(\nu z)}
\cdot \frac{I_{\nu}(\nu z \e)}{K_{\nu}(\nu z \e)} \right)
=\log \left( 1- e^{2\nu(\xi_{\e}-\xi)}(1+\kappa_2(\nu,z))\right),
\end{split}
\end{equation}
We need to consider the difference $(\xi_{\e}-\xi)$ in detail.
\begin{align*}
\xi_{\e}-\xi=\sqrt{1+(\e z)^2}-\sqrt{1+z^2}+\log\left(\frac{\e z}{1+\sqrt{1+(\e z)^2}}\right) - 
\log\left(\frac{ z}{1+\sqrt{1+z^2}}\right)&\\
=\sqrt{1+(\e z)^2}\left[1-\frac{1}{\e}\sqrt{\frac{\e^2+(\e z)^2}{1+(\e z)^2}}\right] +
\log\left(\frac{\e z}{1+\sqrt{1+(\e z)^2}}\right) - \log\left(\frac{ z}{1+\sqrt{1+z^2}}\right)&.
\end{align*}
We are interested in the asymptotic behaviour of $(\xi_{\e}-\xi)$ as $\e \to 0$, which is possibly
non-uniform in $z$. Hence, we consider $(\xi_{\e}-\xi)$ under three asymptotic regimes, 
$|\e z|\to \infty, |\e z|\to 0$ and $|\e z|\sim \textup{const}$.  We find by straightforward estimates
\begin{equation}\begin{split}
&\textup{Re}\,(\xi_{\e}-\xi) \sim \e \, \textup{Re}(z)(1-1/\e) = 
\textup{Re}(z)(\e-1), \ \textup{as} \ |\e z|\to \infty, \ \e \to 0,\\
&\textup{Re}\,(\xi_{\e}-\xi) \sim \log |\e z|-\textup{Re}\sqrt{1+z^2}, \ \textup{as} \ |\e z|\to 0, \ \e \to 0,\\
&\textup{Re}\,(\xi_{\e}-\xi) \sim -C\e^{-1}, \ \textup{as} \ |\e z|\sim \textup{const}, \ \e \to 0,
\end{split}\end{equation}
for some constant $C>0$. For $\{z\in \C||\textup{arg}(z)|<\pi /2\}\cup \{z=ix|x\in (-1,1)\}$, 
we have $\textup{Re}\sqrt{1+z^2}>0$, and $\textup{Re}(z)>0$ as $|z|\to \infty$. 
Consequently, for $\e>0$ sufficiently small
$\textup{Re}(\xi_{\e}-\xi)<\delta <0$
for some fixed $\delta<0$ and hence $\exp(2\nu(\xi_{\e}-\xi))$ 
vanishes as $\nu\to \infty$, uniformly in any compact subset of 
$\{z\in \C||\textup{arg}(z)|<\pi /2\}\cup \{ix|x\in (-1,1)\}$. 

Note also
\begin{equation}\label{log4}\begin{split}
&\log \left(1-\frac{K_{\nu}(\nu z)}{I_{\nu}(\nu z)}\cdot \frac{ \nu z \e I'_{\nu}(\nu z\e)
+ \A_k I_{\nu}(\nu z \e)}{ \nu z\e K'_{\nu}(\nu z\e) + \A_k K_{\nu}(\nu z\e)}\right)
= \log \left( 1- e^{-2\nu(\xi_{\e}+\xi)}(1+\kappa_3(\nu,z))\right),\\
&\log \left(1-\frac{K_{\nu}(\nu z)}{I_{\nu}(\nu z)}\cdot \frac{ \nu z \e I'_{\nu}(\nu z\e)
- \A_k I_{\nu}(\nu z \e)}{ \nu z\e K'_{\nu}(\nu z\e) - \A_k K_{\nu}(\nu z\e)}\right)
= \log \left( 1- e^{-2\nu(\xi_{\e}+\xi)}(1+\kappa_4(\nu,z))\right), 
\end{split}\end{equation}
where the error terms $\kappa_3(\nu,z)$ and $\kappa_4(\nu,z)$ are bounded for large $\nu$ 
uniformly in any compact subset of $\{z\in \C||\textup{arg}(z)|<\pi /2\}\cup \{ix|x\in (-1,1)\}$;
and $\textup{Re}(\xi_{\e}+\xi)>0$.
The uniform expansions above show 
that in \eqref{log1}, \eqref{log2}, \eqref{log3} and \eqref{log4} the arguments of the logarithms stay away 
from the branch cut $\C\backslash \R^-$ for $\nu$ large enough and $\e>0$ 
sufficiently small, uniformly in any compact subset of 
$\{z\in \C||\textup{arg}(z)|<\pi /2\}\cup \{ix|x\in (-1,1)\}$. Consequently, in view of the expression
\eqref{t-bessel}, $t_{\eta, \e}^{k}(\lambda)$ is in particular holomorphic 
in an open neighborhod of $\{\lambda \in [0,c']\}\subset \C$
for some $0\!<\!c'\!<\!1$ and $\nu(\eta)>\nu_0$. 
For any $\eta\in E_k$, $t_{\eta, \e}^{k}(\lambda)$ is moreover holomorphic 
in $\lambda \in \C\backslash \{x\in \R\mid x>c(\eta)\}$. Thus, setting 
$c:=\min \{c', c(\eta)\mid \eta \in E_k, \nu(\eta)
\leq \nu_0\},$ we deduce for $\e>0$ 
sufficiently small
\begin{align*}
\zeta_k(s,\e)=\sum_{\eta \in E_k}  \nu(\eta)^{-2s} 
\frac{s^2}{\Gamma(s+1)}\int_0^{\infty}\frac{t^{s-1}}{2\pi i}
\int_{\wedge_{c}}\frac{e^{-\lambda t}}{-\lambda}\, 
 t_{\eta, \e}^{k}(\lambda)  d\lambda \, dt,
\end{align*}
where deforming the integration contour from $\Lambda_{c(\eta)}$ to $\Lambda_c$ 
is permissible, as the deformation in performed within the region of regularity for each 
$t_{\eta, \e}^{k}(\lambda), \eta\in E_k$. Employing again the expansions 
\eqref{large-nu-I} and \eqref{large-nu-K} we find that 
\begin{align}
\sum_{\eta \in E_k}  t_{\eta, \e}^{k}(\lambda) \nu(\eta)^{-2s}, \ 
\textup{Re}(s)\gg 0,
\end{align}
converges uniformly in $\lambda \in \Lambda_c$ and hence by the 
uniform convergence of the integrals and series we arrive 
at the statement of the proposition.
\end{proof}

\begin{prop}\label{large-nu} 
Consider the notation fixed in Proposition \ref{N-prop} and \ref{contour-change}. 
Let $\lambda \in \Lambda_c$ 
and $t_{\e}(\lambda):=(1-(\e^2 \lambda))^{-1/2}$. 
Then for $\e>0$ sufficiently small we have the following 
asymptotic expansion for large $\nu(\eta), \eta\in E_k$
\begin{align*}
t_{\nu,\e}^{k}(\lambda) \sim  \sum_{r=1}^{\infty}(-\nu)^{-r}\left(
-2D_r(t_{\e})+ M_r(t_{\e},-\A_k)+M_r(t_{\e},\A_k) - \frac{(-1)^{r+1}}{r}(\A_k^r + (-\A_k)^r)\right).
\end{align*}
\end{prop}

\begin{proof}
We discuss the expansions of the individual terms in the expression for 
$t_{\nu,\e}^{k}(\lambda)$ in Proposition \ref{N-prop}. Using the expansions 
\eqref{large-nu-I}, \eqref{large-nu-K} and \eqref{polynom2} we compute for 
large $\nu\in F_k$, putting $\eta_{\e}:=1/t_{\e}+\log (\e z/(1+1/t_{\e}))$
\begin{align*}
\log \left(\frac{\pm \A_k}{\nu} K_{\nu}(\nu z \e) -z \e K'_{\nu}(\nu z \e)\right) 
&\sim \log \sqrt{\frac{\pi}{2\nu}}\frac{e^{-\nu \eta_{\e}}}{(1+z^2\e^2)^{-1/4}} + 
\sum_{r=1}^{\infty}\frac{M_r(t, \pm \A_k)}{(-\nu)^r}, \\
\log K_{\nu}(\nu z \e)
&\sim \log \sqrt{\frac{\pi}{2\nu}}\frac{e^{-\nu \eta_{\e}}}{(1+z^2\e^2)^{1/4}} + 
\sum_{r=1}^{\infty}\frac{D_r(t)}{(-\nu)^r}.
\end{align*}

The standard expansion of the logarithm yields
\begin{align*}
\log\left(1-\frac{\A^2_k}{\nu^2}\right)=\log\left(1+\frac{\A_k}{\nu}\right)+
\log\left(1-\frac{\A_k}{\nu}\right)=\sum_{r=1}^{\infty}(-1)^{r+1}\frac{(\A_k^r+(-\A_k)^r)}{r(-\nu)^r}.
\end{align*}

This already gives all the terms in the stated asymptotic expansion of $t_{\nu,\e}^{k}(\lambda)$. 
Thus we need to check that the remaining terms indeed have no asymptotic contribution. 
The remaining terms are estimated, using \eqref{large-nu-I}, \eqref{large-nu-K}, putting 
$\xi_{\e}:=1/t_{\e}(\lambda)+\log (\e z/(1+1/t_{\e}(\lambda)))$, as follows
\begin{align}\label{eta-difference}
&\frac{ K_{\nu}(\nu z)}{ I_{\nu}(\nu z)} \cdot \frac{ \nu z \e I'_{\nu}
(\nu z \e)\pm\A_k I_{\nu}(\nu z \e)}{ \nu z \e K'_{\nu}(\nu z \e)\pm
\A_k K_{\nu}(\nu z \e)}\sim_{\nu\to\infty} O(e^{2\nu(\eta_{\e}-\eta)}), \\
&\frac{ \nu z K'_{\nu}(\nu z)\pm\A_k K_{\nu}(\nu z)}{ \nu z I'_{\nu}(\nu z)
\pm\A_k I_{\nu}(\nu z)}\cdot \frac{ I_{\nu}(\nu z \e)}{ K_{\nu}(\nu z \e)}
\sim_{\nu\to\infty} O(e^{2\nu(\eta_{\e}-\eta)}).
\end{align}

The difference $(\xi_{\e}-\xi)$ has been considered in detail
in Proposition \ref{contour-change}.
For $\e$ sufficiently small, $\textup{Re}(\xi_{\e}-\xi)<0$ and hence $O(e^{2\nu(\xi_{\e}-\xi)})$ 
in \eqref{eta-difference} does not contribute to the asymptotic expansion for large $\nu$.
\end{proof}

Let us introduce a new (shifted) zeta-function
\begin{align}\label{shifted-zeta}
\zeta_{k,N}(s):=\sum_{\eta \in E_k}\nu(\eta)^{-s}=\zeta\left(\frac{s}{2},\, \Delta_{k,\ccl,N}+\A_k^2\right), \, \textup{Re}(s)>n.
\end{align}
The heat trace expansions for $(\Delta_{k,\ccl,N}+\A_k^2)$ and $\Delta_{k,\ccl,N}$ have 
the same exponents, and hence $\zeta_{k,N}(s)$ extends meromorphically to $\C$ with simple poles at 
$\{(n-2k)\mid k\in \N\}$. Consequently, terms $\nu(\eta)^{-r}$ in the asymptotic expansion of $t_{\eta,\e}^{k}(\lambda)$ 
with $r=n-2k,k\in \N$, may lead to singular behaviour of $T^k_{\e}(s,\lambda)$ at $s=0$. 
In particular the poles occur only at the odd integers, since $\dim N$ is odd. 
We regularize $T^k_{\e}(s,\lambda)$ by subtracting off these terms from $t_{\eta,\e}^{k}(\lambda)$, and define
\begin{equation}\label{P-k}
\begin{split}
f_{r,\e}^{k}(\lambda):=& 2D_{2r+1}(t_{\e})- M_{2r+1}(t_{\e},-\A_k)-M_{2r+1}(t_{\e},\A_k), \\
p_{\eta,\e}^{k}(\lambda):=&t_{\eta,\e}^{k}(\lambda)-\sum_{r=1}^{(n-1)/2}\nu(\eta)^{-(2r+1)} 
f_{r,\e}^{k}(\lambda), \quad P^k_{\e}(s,\lambda):=\sum_{\eta\in E_k}p_{\eta,\e}^{k}(\lambda)
\nu(\eta)^{-2s}.
\end{split}
\end{equation}
By construction, $P^k_{\e}(s,\lambda)$ is regular at $s=0$. Contribution of the terms $f_{r,\e}^{k}(\lambda)$ 
is computed in terms of the polynomial structure of the coefficients $M_r(t,\A)$ and $D_r(t)$ in \eqref{MD-polynom}. 
The computation uses special integrals evaluated already by Spreafico \cite{Spr:ZIF}.

\begin{prop}\label{ff}
\begin{align*}
\int_0^{\infty}t^{s-1}\frac{1}{2\pi i}\int_{\wedge_c}\frac{e^{-\lambda t}}
{-\lambda}\, f_{r,\e}^{k}(\lambda)\, d\lambda \, dt =
\sum_{b=0}^{2r+1}\left(2x_{2r+1,b} - z_{2r+1,b}(-\A_k)-z_{2r+1,b}(\A_k)\right) \\ \times 
\frac{\Gamma(s+b+r+1/2)}{s\Gamma (b+r+1/2)}\e^{2s}.
\end{align*}
\end{prop}

\begin{proof} 
The $\lambda$-independent part of $f_{r,\e}^{k}(\lambda)$ vanishes after integration in $\lambda$. 
The coefficients $M_{r}(t_{\e}(\lambda),\pm\A_k)$ and $D_r(t)$ in the definition of $f_{r,\e}^{k}(\lambda)$ are polynomial 
in $t_{\e}(\lambda)=(1-\e^2\lambda)^{-1/2}$. Repeating (\!\!\cite{MV}, Proposition 5.9) or \cite{Spr:ZIF}
we find
\begin{align*}
\int_0^{\infty}t^{s-1}\frac{1}{2\pi i}\int_{\wedge_c}\frac{e^{-\lambda t}}
{-\lambda}\frac{1}{(1-\e^2\lambda)^a}\, d\lambda \, dt 
= \e^{2s} \frac{\Gamma(s+a)}{s\, \Gamma (a)}.
\end{align*}
Polynomial representation \eqref{MD-polynom} yields the statement.
\end{proof}

Consequently we arrive at the intermediate 
representation of $\zeta_k(s,\e)$ for $\textup{Re}(s)\gg 0$
\begin{equation}\label{zeta-intermediate}
\begin{split}
\zeta_k(s,\e)&=
\frac{s^2}{\Gamma(s+1)}\int_0^{\infty}\frac{t^{s-1}}{2\pi i}
\int_{\wedge_{c}}\frac{e^{-\lambda t}}{-\lambda}\, 
P^k_{\e}(s,\lambda) d\lambda \, dt \\
&+ \sum_{r=1}^{(n-1)/2}\zeta_{k,N}(2s+2r+1)\frac{s}{\Gamma (s+1)}
\sum_{b=0}^{2r+1} \frac{\Gamma(s+b+r+1/2)}{s\Gamma (b+r+1/2)}\, \e^{2s} \\
& \times \left(2x_{2r+1,b} - z_{2r+1,b}(-\A_k)-z_{2r+1,b}(\A_k)\right).
\end{split}
\end{equation}

While the second summand in \eqref{zeta-intermediate} 
extends meromorphically to $\C$, 
it still remains to derive an analytic extension to $s=0$ 
for the first summand.

\begin{prop}\label{AB} 
Consider notation fixed in Proposition \ref{N-prop} and \eqref{P-k}. 
Then for large arguments $\lambda\in \Lambda_c$ and fixed order $\nu$ we have the following asymptotics
\begin{align*}
p_{\nu,\e}^{k}(\lambda)=a^k_{\nu,\e}\log (-\lambda)+b^k_{\nu,\e}+O\left((-\lambda)^{-1/2}\right),
\end{align*}
where 
$$
a^k_{\nu}=1, \quad b^k_{\nu}=2\log \e - \log \left(1-\frac{\A_k^2}{\nu^2}\right).
$$ 
\end{prop}
\begin{proof}
The function $p_{\eta,\e}^{k}(\lambda)$ is given by the following expression
\begin{align*}
p_{\eta,\e}^{k}(\lambda)=t_{\eta,\e}^{k}(\lambda)-\sum_{r=1}^{(n-1)/2}\nu(\eta)^{-(2r+1)} 
2D_{2r+1}(t_{\e})- M_{2r+1}(t_{\e},-\A_k)-M_{2r+1}(t_{\e},\A_k).
\end{align*}
The polynomials $M_{2r+1}(t_{\e}(\lambda),\pm \A_k)$ and $D_{2r+1}(t)$ have no constant terms, and hence 
are $O\left((-\lambda)^{-1/2}\right), \lambda \to \infty,$ since
\begin{align}
t_{\e}(\lambda)=\frac{1}{\sqrt{1-\e^2\lambda}}=O\left((-\lambda)^{-1/2}\right), 
\quad  \lambda \to \infty.
\end{align} 
By \eqref{large-arg-I} and \eqref{large-arg-K}, setting $\nu\equiv \nu(\eta)$ we find as $\lambda \to \infty$
\begin{align*}
\log \left(\pm \frac{\A_k}{\nu} K_{\nu}(\nu z \e)-z\e K'_{\nu}(\nu z \e)\right)
&\sim \log \sqrt{\frac{\pi}{2z\nu \e}} + \log z\e + \log \left(1\pm \frac{\A_k}{\nu z \e}\right), \\
\log \left(K_{\nu}(\nu z \e)\right) &\sim \log \sqrt{\frac{\pi}{2z\nu \e}}
\end{align*}

Moreover we have
\begin{align*}
&\frac{ \nu z K'_{\nu}(\nu z)\pm\A_k K_{\nu}(\nu z)}{ \nu z I'_{\nu}(\nu z)\pm\A_k 
I_{\nu}(\nu z)}\cdot \frac{ I_{\nu}(\nu z \e)}{ K_{\nu}(\nu z \e)}\sim_{\lambda \to\infty} O(e^{2\nu z(\e-1)}), \\
&\frac{ K_{\nu}(\nu z)}{ I_{\nu}(\nu z)} \cdot \frac{ \nu z \e I'_{\nu}(\nu z \e)\pm\A_k 
I_{\nu}(\nu z \e)}{ \nu z \e K'_{\nu}(\nu z \e)\pm\A_k K_{\nu}(\nu z \e)}\sim_{\lambda \to\infty} O(e^{2\nu z(\e-1)}).
\end{align*}

$(\e-1)<0$ and Re$(z)>0$ for large $z=\sqrt{-\lambda}, \lambda \in \Lambda_c$. Consequently 
$O(e^{2\nu z(\e-1)})$ is in particular of $O\left((-\lambda)^{-1/2}\right)$ asymptotics 
for $\lambda \to \infty, \lambda \in \Lambda_c$. By the explicit expression for $t_{\eta,\e}^{k}(\lambda)$ 
in \eqref{t-bessel} the statement follows.
\end{proof}

\begin{defn} \label{AB1} Define for $\textup{Re}(s)>n$ in notation of Proposition \ref{AB}
\begin{align}
A^k_{\e}(s):=\sum_{\eta \in E_k}a^k_{\eta,\e}\, \nu^{-2s}, \quad B^k_{\e}(s):=\sum_{\eta \in E_k}b^k_{\eta,\e}\, \nu^{-2s}.
\end{align}
\end{defn}

\begin{prop}\label{P}
Consider notation fixed in Proposition \ref{N-prop} and \eqref{P-k}. Then
\begin{align*}
P^k_{\e}(s,0)=0.
\end{align*}
\end{prop}
\begin{proof}
By \eqref{DM}
\begin{align}
M_{r}(1,-\A_k)- M_{r}(1,\A_k)=\frac{(-\A_k)^r-\A_k^r}{r}.
\end{align}
For any fixed $\e>0$ clearly $\lambda \to 0$ implies that $t=(1-\e^2\lambda)^{-1/2}$ tends to $1$. 
Hence 
\begin{align*}
f_{r,\e}^{k}(0)=(-1)^{2r+1}\frac{(-\A_k)^{2r+1}+\A_k^{2r+1}}{r}=0.
\end{align*}

Moreover, by \eqref{small}
\begin{align*}
\log \left(\pm \frac{\A_k}{\nu} K_{\nu}(\nu z \e)-z\e K'_{\nu}(\nu z \e)\right)
&\sim \log [2^{\nu-1}\Gamma(\nu)(\nu z\e)^{-\nu}] + \log \left(1\pm \frac{\A_k}{\nu}\right), \\
\log \left(K_{\nu}(\nu z \e)\right) &\sim \log [2^{\nu-1}\Gamma(\nu)(\nu z\e)^{-\nu}].
\end{align*}

Moreover we have
\begin{align*}
&\frac{ \nu z K'_{\nu}(\nu z)\pm\A_k K_{\nu}(\nu z)}{ \nu z I'_{\nu}(\nu z)\pm
\A_k I_{\nu}(\nu z)}\cdot \frac{ I_{\nu}(\nu z \e)}{ K_{\nu}(\nu z \e)}\sim_{\lambda \to 0} \e^{2\nu}, \\
&\frac{ K_{\nu}(\nu z)}{ I_{\nu}(\nu z)} \cdot \frac{ \nu z \e I'_{\nu}(\nu z \e)
\pm\A_k I_{\nu}(\nu z \e)}{ \nu z \e K'_{\nu}(\nu z \e)\pm\A_k K_{\nu}(\nu z \e)}
\sim_{\lambda \to 0} \e^{2\nu}.
\end{align*}

By the explicit expression for $t_{\eta,\e}^{k}(\lambda)$ in \eqref{t-bessel} the statement 
follows. Note that the $\e-$dependence cancels.
\end{proof}

We have now all ingredients together to write down the meromorphic continuation to $s=0$
of the zeta-function $\zeta_k(s,\e)$, introduced in Proposition \ref{N-prop}. 
By the arguments of (\!\!\cite{S}, Section 4.1) we have
\begin{equation}\label{zeta-expression}
\begin{split}
\zeta_k(s,\e) &=  \frac{s}{\Gamma (s+1)}[\gamma A^k_{\e}(s)-B^k_{\e}(s)-\frac{1}{s}A^k_{\e}(s)+P^k_{\e}(s,0)] \\ 
&+  \sum_{r=1}^{(n-1)/2}\frac{s^2}{\Gamma (s+1)}\zeta_{k,N}(2s+2r+1)\int_0^{\infty}t^{s-1}\frac{1}{2\pi i}
\int_{\wedge_c}\frac{e^{-\lambda t}}{-\lambda}f^k_{r,\e}(\lambda)\, d\lambda \, dt \\
&+ \frac{s^2}{\Gamma (s+1)}h(s,\e),
\end{split}
\end{equation}
where $h(s,\e)$ vanishes with its derivative at $s=0$. Note that all the terms are regular
at $s=0$. Inserting the results of Proposition \ref{ff}, Proposition \ref{AB}, 
Proposition \ref{P} together with Definition \ref{AB1} into the expression \eqref{zeta-expression} we obtain the following

\begin{prop}\label{zeta-total} 
Let $(E_N,\nabla_N,h_N)$ be a flat Hermitian vector bundle over an even-dimensional 
oriented closed Riemannian manifold $(N^n,g^N)$.
Denote by $\Delta_{k,\ccl,N}$ the Laplacian on coclosed $k-$differential forms 
$\Omega^k_{\textup{\ccl}}(N,E_N)$. Consider notation fixed in \eqref{MD-polynom} and \eqref{shifted-zeta}. Put
\begin{align*}
&\A_k:=\frac{(n-1)}{2}-k, \quad 
F_k:=\left\{\nu=\sqrt{\eta + \A_k^2} \mid \eta \in \textup{Spec}\Delta_{k,\ccl,N}\backslash \{0\}\right\}, \\
&\zeta_{k,N}(s):=\sum_{\nu \in F_k} \nu^{-s}, \quad Re(s) \gg 0.
\end{align*}
Then for $\e>0$ sufficiently small, $\zeta_k(s,\e)$ defined in Definition \ref{zetas}
admits an analytic continuation to $s=0$ of the form
\begin{equation}\label{zeta-continued}
\begin{split}
\zeta_k(s,\e)=\frac{s}{\Gamma (s+1)}\left[ \left((\gamma -2 \log \e  -\frac{1}{s}\right)\zeta_{k,N}(2s)  + 
\sum_{\nu \in F_k}\nu^{-2s} \log\left(1-\frac{\A_k^2}{\nu^2}\right) \right] \\
+\sum_{r=1}^{(n-1)/2}\zeta_{k,N}(2s+2r+1)\frac{s}{\Gamma (s+1)}
\left(\sum_{b=0}^{2r+1}\left[2x_{2r+1,b}- z_{2r+1,b}(-\A_k)-z_{2r+1,b}(\A_k)\right]\right. \\
\left. \times \frac{\Gamma(s+b+r+1/2)}{\Gamma (b+r+1/2)}\right)\e^{2s} +s^2h(s) / \Gamma (s+1),
\end{split}
\end{equation}
where $h(s)$ vanishes with its derivative at $s=0$.
\end{prop}

Note the full analogy (up to computationally irrelevant, but geometrically crucial sign differences) 
to the corresponding result in (\!\!\cite{Ver:ATO}, Proposition 6.10). An ad verbatim repetition of the 
arguments in the proof of (\!\!\cite{Ver:ATO}, Corollary 6.11) leads to the final formula.
\begin{cor}\label{total-contribution-1}
Let $(E_N,\nabla_N,h_N)$ be a flat Hermitian vector bundle over an even-dimensional 
oriented closed Riemannian manifold $(N^n,g^N)$.
Denote by $\Delta_{k,\ccl,N}$ the Laplacian on coclosed $k-$differential forms 
$\Omega^k_{\textup{\ccl}}(N,E_N)$ and put
\begin{align*}
&\A_k:=\frac{(n-1)}{2}-k, \quad 
F_k:=\left\{\nu=\sqrt{\eta + \A_k^2} \mid \eta \in \textup{Spec}
\Delta_{k,\ccl,N}\backslash \{0\}\equiv E_k\right\}, \\
&\zeta_{k,N}(s)=\sum_{\nu \in F_k} \nu^{-s},\quad \zeta(s, \Delta_{k,ccl,N}):=
\sum_{\nu\in E_k} \eta^{-s}, \quad Re(s)\gg0.
\end{align*}
Then we find in notation of \eqref{MD-polynom} for $\e>0$ sufficiently small
\begin{align*}
\zeta_k'(0,\e)&= - \zeta' (0, \Delta_{k,\ccl,N}) -2 \log \e \cdot \zeta (0, \Delta_{k,\ccl,N}) \\
&+\frac{1}{2}\sum_{r=1}^{(n-1)/2}\textup{Res}\, \zeta_{k,N}(2r+1)\left(\sum_{b=0}^{2r+1}
\left[2x_{2r+1,b}- z_{2r+1,b}(-\A_k)-z_{2r+1,b}(\A_k)\right]\right. \\ &\left. \times 
\frac{\Gamma'(b+r+1/2)}{\Gamma (b+r+1/2)}\right). 
\end{align*}
\end{cor}

\begin{proof}
We begin with the following observation
\begin{align*}
\sum_{\nu \in F_k}\nu^{-2s} \log\left(1-\frac{\A_k^2}{\nu^2}\right) &= 
\sum_{\nu \in F_k}\nu^{-2s} \log\left(\frac{\nu^2-\A_k^2}{\nu^2}\right)\\
&= 2 \, \zeta'_{k,N}(2s) + \sum_{\eta \in E_k}(\eta + \A_k^2)^{-s} \log \eta \\
&= 2 \, \zeta'_{k,N}(2s)  - \sum_{j=0}^{\infty} {-s \choose j} \A_k^{2j} \zeta' (s+j, \Delta_{k,\ccl,N})
\end{align*}
Since the poles of $\zeta (s, \Delta_{k,\ccl,N})$ lie on half integers for $n=\dim N$ odd, 
the summands in the sum above are regular at $s=0$ and hence we find
\begin{align*}
\left. \frac{d}{ds}\right|_{s=0} \frac{s}{\Gamma (s+1)}
\sum_{\nu \in F_k}\nu^{-2s} \log\left(1-\frac{\A_k^2}{\nu^2}\right) = 
- \zeta' (0, \Delta_{k,\ccl,N}) + 2 \, \zeta'_{k,N}(0). 
\end{align*}
Similarly we find
\begin{align*}
\zeta_{k,N}(2s) = 
\sum_{\nu \in F_k}\nu^{-2s} = \sum_{\eta \in E_k}(\eta + \A_k^2)^{-s}
= \sum_{j=0}^{\infty} {-s \choose j} \A_k^{2j} \zeta (s+j, \Delta_{k,\ccl,N}). 
\end{align*}
Again, by regularity of the summands at $s=0$, we deduce
$$\zeta_{k,N}(0)= \zeta (0, \Delta_{k,\ccl,N}).$$
Moreover we compute
\begin{align*}
\left. \frac{d}{ds}\right|_{s=0} \frac{s}{\Gamma (s+1)}\zeta_{k,N}(2s) \left(\gamma -\frac{1}{s}\right) = - 2 \zeta'_{k,N}(0). 
\end{align*}
The statement now follows by a combination of the two observations
\begin{align*}
\left.\frac{d}{ds}\right|_{s=0}\zeta_{k,N}(2s+i)\frac{s}{\Gamma (s+1)}\frac{\Gamma(s+b+i/2)}{\Gamma (b+i/2)}=\\=\frac{1}{2}\,\textup{Res}\zeta_{k,N}(i)\left[\frac{\Gamma'(b+i/2)}{\Gamma (b+i/2)}+\gamma \right]+\textup{PP}\zeta_{k,N}(i), 
\end{align*}
where PP$\zeta_{k,N}(i)$ denotes the constant term in the asymptotics of $\zeta_{k,N}(s)$ near the pole singularity $s=i$. The second observation is
\begin{align*}
\sum_{b=0}^{2r+1}&\left(2x_{2r+1,b}-z_{2r+1,b}(-\A_k)-z_{2r+1,b}(\A_k)\right)\\
&= 2D_{2r+1}(1)-M_{2r+1}(1,-\A_k) - M_{2r+1}(1,\A_k)\\
&=(-1)^{2r+1}\frac{\A_k^{2r+1}+(-\A_k)^{2r+1}}{2r+1}=0.
\end{align*}
\end{proof}

\begin{prop}\label{h-truncated}
Let $(E_N,\nabla_N,h_N)$ be a flat Hermitian vector bundle over an even-dimensional 
oriented closed Riemannian manifold $(N^n,g^N)$. Denote the Euler characteristic of 
$(N,E_N)$ by $\chi(N,E_N)$ and the Betti numbers by $b_k:=\dim H^k(N,E_N)$. Then 
in notation of Definition \ref{zetas} we find
\begin{equation}
\begin{split}
\sum_{k=0}^{n}\frac{(-1)^{k+1}}{2}\, \zeta'_{k,H}(0,\e)&=
\frac{1}{2}\log \e \sum_{k=0}^n (-1)^{k} \, k \, b_k
- \frac{1}{2}\sum_{k=0}^{(n-1)/2}(-1)^kb_k \log (n-2k+1).
\end{split}
\end{equation}
\end{prop}

\begin{proof}
By Definition \ref{zetas} we can write
\begin{equation}
\begin{split}
\sum_{k=0}^{n}\frac{(-1)^{k+1}}{2}\, \zeta_{k,H}(s,\e)&=
\sum_{k=0}^n\frac{(-1)^{k+1}}{2} \, b_k \, \zeta(s, H^k_{0,\e,\textup{\textup{mix}}})\\ 
&-\sum_{k=0}^n\frac{(-1)^{k+1}}{2} \, b_k \, \zeta(s, H^k_{0,\textup{\textup{mix}}})=:H(s,\e)-H(s).
\end{split}
\end{equation}
$H'(0)$ has been evaluated in (\!\!\cite{Ver:ATO}, Theorem 7.8) with 
\begin{equation}\label{H1}
\begin{split}
\sum_{k=0}^n\frac{(-1)^{k+1}}{2}\, b_k \, \zeta'(0, H^k_{0,\textup{\textup{mix}}})=
\frac{1}{2}\sum_{k=0}^{(n-1)/2}(-1)^kb_k \log (n-2k+1).
\end{split}
\end{equation}
We evaluate $H'(0,\e)$ using  (\!\!\cite{Les:DOR}, Theorem 1.2), which relates the zeta determinants 
to the normalized solutions of the operators, satisfying the corresponding boundary conditions. 
The boundary conditions for $H^k_{0,\e,\textup{\textup{mix}}}$ 
have been determined in Proposition \ref{rel-bc-prop}. 
By the formula (\!\!\cite{Les:DOR}, Theorem 1.2) we then find
\begin{align}
\det\nolimits_{\zeta} \left(H^k_{0,\e,\textup{\textup{mix}}}\right) = 
2 \, \e^{k-n/2}.
\end{align}
Taking logarithms and employing Poincare duality on $(N,g^N)$ we find
\begin{equation}\label{H2}
\begin{split}
H'(0,\e) =\log \e \sum_{k=0}^n\frac{(-1)^{k}}{2} \cdot k \cdot b_k.
\end{split}
\end{equation}
The statement follows by combination of \eqref{H1} and \eqref{H2}.
\end{proof}

Summing up the expressions in Corollary \ref{total-contribution-1} and Proposition \ref{h-truncated}, 
we arrive at the following result.
\begin{thm}\label{t-difference}
Let $(\mathscr{C}(N)=(0,1]\times N, g=dx^2\oplus x^2g^N)$ be an odd-dimensional 
bounded cone over a closed oriented Riemannian manifold $(N^n,g^N)$. 
Denote by $(\mathscr{C}_\e (N)=[\e,1]\times N, g)$ 
its truncation. Let $(E,\nabla, h^E)$ be a flat complex Hermitian vector bundle 
over $(\mathscr{C}(N),g)$ and $(E_N,\nabla_N,h_N)$ the corresponding restriction 
to the cross-section $N$. Denote by $\chi(N,E_N)$ the Euler characteristic and by $b_k:=\dim H^k(N,E_N)$ 
the Betti numbers of $(N,E_N)$. Denote by $\Delta_{k,\ccl,N}$ the Laplacian on coclosed $k-$differential forms 
$\Omega^k_{\textup{\ccl}}(N,E_N)$ and put
\begin{align*}
&\A_k:=\frac{(n-1)}{2}-k, \quad 
F_k:=\left\{\nu=\sqrt{\eta + \A_k^2} \mid \eta \in \textup{Spec}
\Delta_{k,\ccl,N}\backslash \{0\}\right\}, \\
&\zeta_{k,N}(s)=\sum_{\nu \in F_k} \nu^{-s},\quad  
\delta_k:=\left\{ \begin{array}{cl} 1/2 & \textup{if} \ k=(n-1)/2, \\ 1 & \textup{otherwise}.
\end{array}\right.
\end{align*}
Then the difference of the scalar analytic torsions for $(\mathscr{C}(N),g)$ and $(\mathscr{C}(N)_{\e}, g)$ 
is given by the following explicit expression
\begin{align*}
&\log T(\mathscr{C}(N)_{\e}, E, g)-\log T(\mathscr{C}(N), E, g)= 
\frac{1}{2}\log T(N,E_N,g^N) \\ -&\, \frac{1}{2}\sum_{k=0}^{(n-1)/2}(-1)^kb_k \log (n-2k+1) 
+ \sum_{k=0}^{(n-1)/2}\frac{(-1)^k}{4}\delta_k \sum_{r=1}^{(n-1)/2}\textup{Res}\, \zeta_{k,N}(2r+1) \\
\times &\, \sum_{b=0}^{2r+1}\left[2x_{2r+1,b}- z_{2r+1,b}(-\A_k)-z_{2r+1,b}(\A_k)\right] 
\frac{\Gamma'(b+r+1/2)}{\Gamma (b+r+1/2)}. 
\end{align*}
\end{thm}

\begin{proof}
The result is a consequence of Corollary \ref{total-contribution-1}, Proposition \ref{h-truncated} and 
the following two combinatorial identities, which follow from Poincare duality on $N$. First
\begin{align*}
- \sum_{k=0}^{(n-1)/2}\frac{(-1)^k}{2}\delta_k \zeta'(0, \Delta_{k,\ccl,N}) &= \frac{1}{4} 
\sum_{k=0}^n (-1)^{k+1} \zeta'(0, \Delta_{k,\ccl,N}) \\ &=  \frac{1}{4}\sum_{k=0}^n (-1)^k 
k \zeta'(0, \Delta_{k,N})=  \frac{1}{2}\log T(N,E_N,g^N).
\end{align*}
Similarly we obtain
\begin{align*}
- 2\log \e \sum_{k=0}^{(n-1)/2}\frac{(-1)^k}{2}\delta_k \zeta(0, \Delta_{k,\ccl,N}) = \log \e  \sum_{k=0}^n \frac{(-1)^{k+1}}{2} \zeta(0, \Delta_{k,\ccl,N}) \\ =  \log \e \sum_{k=0}^n \frac{(-1)^k}{2} k \zeta(0, \Delta_{k,N})=- \log \e \sum_{k=0}^n \frac{(-1)^k}{2} \cdot k \cdot b_k.
\end{align*}
Comparing this to Proposition \ref{h-truncated} in particular shows that the $\e$-dependence cancels in the total formula.
\end{proof}

Comparison of Theorem \ref{t-difference} and Theorem \ref{BV-Theorem} yields the following
\begin{cor}\label{t1}
\begin{align*}
\log T(\mathscr{C}(N)_{\e}, E, g)&= \sum_{k=0}^{(n-1)/2}\frac{(-1)^k}{2}
\delta_k \sum_{r=1}^{(n-1)/2}\textup{Res}\, \zeta_{k,N}(2r+1) \\
&\times \sum_{b=0}^{2r+1}\left[2x_{2r+1,b}- z_{2r+1,b}(-\A_k)-z_{2r+1,b}(\A_k)\right] 
\frac{\Gamma'(b+r+1/2)}{\Gamma (b+r+1/2)}. 
\end{align*}
\end{cor}

\section{Metric Anomaly at the Regular Boundary of the Cone}
\label{section-anomaly}

Alternatively, the analytic torsion of the cone-like cylinder $(C_{\e}(N),g=dx^2\oplus x^2g^N)$ 
can be computed by relating it to the analytic torsion of the exact cylinder, using the anomaly formula of 
Br\"uning-Ma, see (\!\!\cite{BruMa:AAF2}, Theorem 3.4). 
Introducing new coordinates $y=\log (1/x)$ near the right boundary component $\{x=1\}\times N$ of $C_{\e}(N)$, 
and $z=\log (x/\e)$ near the left boundary component $\{x=\e\}\times N$, 
we can write for $\delta >0$ small the Riemannian metric $g$ as follows
\begin{equation}
\begin{split}
g=e^{-2y}\left(dy^2+g^N\right), \ y \in [0,\delta), \ \textup{near $\{x=1\}\times N$ of $C_{\e}(N)$}, \\
g=\e^2e^{2z}\left(dz^2+g^N\right), \ z \in [0,\delta), \ \textup{near $\{x=\e\}\times N$ of $C_{\e}(N)$}. 
\end{split}
\end{equation}
The Levi-Civita connection $\nabla^{LC}$, 
induced by $g$, defines secondary classes $B_{\e}(\nabla^{LC})$ and $B_{1}(\nabla^{LC})$ 
at the left  $\{x=\e\}\times N$ and the right  $\{x=1\}\times N$ boundary components of $C_{\e}(N)$, respectively.
By Proposition \ref{Scaling} and in view of the explicit formulae in \eqref{RS3}, we deduce
\begin{align}\label{2B}
B(g^N) := B_{1}(\nabla^{LC}) = - B_{\e}(\nabla^{LC}),
\end{align}
independent of $\e>0$. Using (\!\!\cite{BruMa:AAF2}, Theorem 3.4) and the fact, that the analytic 
torsion of the exact cylinder $([\e,1]\times N, g_0=dx^2\oplus g^N)$ with mixed boundary conditions 
is trivial, we find
\begin{align}\label{t2}
\log T(\mathscr{C}(N)_{\e}, E, g)=\textup{rank}(E) \int_N B(g^N).
\end{align}

\begin{remark}\label{unusual}
Had we taken relative boundary conditions on both boundaries of the cone-like cylinder, 
then the Br\"uning-Ma anomaly contributions from both boundaries would cancel. On the analytic 
side of our computations this phenomena would appear by the fact that the analytic torsion of a 
cone-like cylinder then no longer captures the "singular" contributions in Theorem \ref{BV-Theorem}. 
By that reason we had to combine the relative and absolute boundary conditions.
\end{remark}

Comparison of \eqref{t2} and Corollary \ref{t1} instantly leads in view of 
Theorem \ref{BV-Theorem} to our our first main result in Theorem \ref{main1}. 
\begin{thm}\label{main1+}
Let $\mathscr{C}(N)=(0,1]\times N, g=dx^2\oplus x^2g^N$ be an even-dimensional bounded cone 
over a closed oriented Riemannian manifold $(N^n,g^N), n=\dim N$. Let $(E,\nabla, h^E)$ be a flat complex Hermitian vector bundle 
over $(\mathscr{C}(N),g)$ and $(E_N,\nabla_N,h_N)$ its restriction 
to the cross-section $N$. Denote by $b_k:=\dim H^k(N,E_N)$ 
the Betti numbers of $(N,E_N)$. Then the logarithm of the scalar analytic torsion of $(\mathscr{C}(N),g)$ 
is given by
\begin{align*}
\log T(\mathscr{C}(N),E,g)=& \sum_{k=0}^{(n-1)/2} \frac{(-1)^k}{2} b_k\log (n-2k+1) \\ 
-&\, \frac{1}{2}\log T(N,E_N,g^N) + \frac{\textup{rank}(E)}{2} \int_N B(g^N).
\end{align*}
\end{thm}

\begin{cor}\label{main2+}
Let $(\mathscr{C}(N)=(0,1]\times N, g=dx^2\oplus x^2g^N)$ be an even-dimensional bounded cone over 
a closed oriented Riemannian manifold $(N^n,g^N)$. Let $(E,\nabla, h^E)$ be a flat complex 
Hermitian vector bundle over $M$ and $(E_N,\nabla_N,h_N)$ its restriction to $N$. Let the metric $g_0$ on $M$ 
coincide with $g$ near the singularity at $x=0$ and be product $dx^2\oplus g^N$ in an open neighborhood of the 
boundary $\{1\}\times N$.  Then the Ray-Singer analytic torsion norm of $(\mathscr{C}(N),g_0)$ is given by
\begin{align*}
\log \|\cdot \|^{RS}_{(\mathscr{C}(N),E,g_0)} &= \sum_{k=0}^{(n-1)/2} \frac{(-1)^{k}}{2} \dim H^k(N,E_N) \log (n-2k+1) 
\\ &- \frac{1}{2}\log T(N,E_N,g^N) + \log \|\cdot \|_{\det \mathscr{H}^*(\mathscr{C}(N),E)},
\end{align*}
where $ \|\cdot \|_{\det \mathscr{H}^*(\mathscr{C}(N),E)}$ denotes the norm on the determinant line 
$\det \mathscr{H}^*(\mathscr{C}(N),E)$, 
induced from the $L^2(g,h^E)$-norm on square integrable harmonic forms with relative boundary conditions.
\end{cor}

\begin{proof}
The metric anomaly formula of Br\"uning-Ma in\cite{BruMa:AAF} holds also in case of manifolds 
with isolated conical singularities away from the variation region. Consequently
\begin{align*}
\log \left(\frac{\|\cdot \|^{RS}_{(\mathscr{C}(N),E,g_0)}}{\|\cdot \|_{\det \mathscr{H}^*(\mathscr{C}(N),E)}}\right)
= \log \left(\frac{\|\cdot \|^{RS}_{(\mathscr{C}(N),E,g_0)}}{\|\cdot \|^{RS}_{(\mathscr{C}(N),E,g)}}\right)+
\log T(\mathscr{C}(N),E,g) \\= \log T(\mathscr{C}(N),E,g) - \frac{\textup{rank}(E)}{2} \int_N B(g^N).
\end{align*}
In view of Theorem \ref{main1+} we have proved the statement.
\end{proof}

\bibliographystyle{amsplain}

\begin{thebibliography}{10}

\bibitem{AbrSte:HOM} M. Abramowitz, I.A. Stegun ed. \emph{"Handbook of mathematical functions with formulas, 
graphs, and mathematical tables"}, Reprint of the 1972 edition. Dover Publications, Inc., New York, (1992)

\bibitem{BFK:OTD} D. Burghelea, L. Friedlander, T. Kappeler \emph{"On the determinant of elliptic boundary value problems on a line segment"}, 
Proc. Amer. Math. Soc. 123, 3027-3038 (1995)

\bibitem{BGKE:ZFD} M. Bordag, B. Geyer, K. Kirsten, E. Elizalde \emph{"Zeta function determinant of the Laplace Operator on the D-dimensional ball"}, Comm. Math. Phys. 179 (1996), no. 1, 215-234

\bibitem{BKD:HKA} M. Bordag, K. Kirsten, J.S. Dowker \emph{"Heat-kernels and functional determinants on the cone"}, Comm. Math. Phys. 182 (1996), 371-394

\bibitem{BruMa:AAF} J. Br\"{u}ning, X. Ma \emph{"An anomaly-formula for Ray-Singer metrics on manifolds with boundary"}, Geom. Funct. An. 16 (2006), No. 4, 767-837

\bibitem{BruMa:AAF2} J. Br\"{u}ning, X. Ma \emph{"On the Gluing Formula for The Analytic Torsion"}, preprint on http://www.math.jussieu.fr/~ma/mypubli/glutc.pdf

\bibitem{BruSee:AIT} J. Br\"{u}ning, R. Seeley \emph{"An index theorem for first order regular singular operators"}, Amer. J. Math 110 (1988), 659-714

\bibitem{Che:ATA} J. Cheeger \emph{"Analytic Torsion and Reidemeister Torsion"}, Proc. Nat. Acad. Sci. USA 74 (1977), 2651-2654

\bibitem{Che:SGOlv:AAS} J. Cheeger \emph{"Spectral Geometry of singular Riemannian spaces"}, J. Diff. Geom. 18 (1983), 575-657 

\bibitem{Dar:IRT} A. Dar \emph{"Intersection R-torsion and analytic torsion for pseudo-manifolds"} Math. Z. 194 (1987), 193-216 

\bibitem{Fr} W. Franz \emph{"\"{U}ber die Torsion einer \"{U}berdeckung"}, J. reine angew. Math. 173 (1935), 245-254 

\bibitem{GRJ:TOI} I.S: Gradsteyn, I.M Ryzhik, Alan Jeffrey \emph{"Table of integrals, Series and Products"}, 5th edition, Academic Press, Inc. (1994)

\bibitem{HS} L. Hartmann, M. Spreafico \emph{"The Analytic Torsion of the Cone over an Odd Dimensional Manifold"}, (2010) preprint on arXiv:math.DG/1001.4755v1

\bibitem{HSpr:ZIF} L. Hartmann, M. Spreafico \emph{"The Analytic Torsion of the Cone over an Odd Dimensional Manifold"},
J. Geom. Phys. 61 (2011) 624-657

\bibitem{HS2} L. Hartmann, M. Spreafico \emph{"An extension of the Cheeger-M\"uller theorem for a cone"},
arXiv:1008.2987v1 [math.DG] (2010)

\bibitem{Les:DOR} M. Lesch \emph{"Determinants of regular singular Sturm-Liouville operators"},  Math. Nachr.  194  (1998), 139--170

\bibitem{L:P} M. Lesch \emph{"The analytic torsion of the model cone"}, Columbus University (1994), unpublished notes.

\bibitem{Lu} W. L\"{u}ck \emph{"Analytic and topological torsion for manifolds with boundary and symmetry"}, J. Diff. Geom. 37 (1993), 263-322

\bibitem{MHS} T. de Melo, L. Hartmann, M. Spreafico, \emph{"Reidemeister torsion and analytic torsion of discs"}, Boll. Unione Mat. Ital. (9) 2 (2009), no. 2, 529--533

\bibitem{Mu} W. M\"uller \emph{"Analytic torsion and $R$-torsion of Riemannian manifolds"}, Adv. in Math.  28  (1978), no. 3, 233--305.

\bibitem{MV} W. M\"uller and B. Vertman \emph{"The Metric Anomaly of Analytic Torsion on Manifolds with Conical Singularities"}, arXiv:1004.2067v3 [math.SP] (2010)

\bibitem{Olv:AAS} F.W. Olver \emph{"Asymptotics and special functions"} AKP Classics, A K Peters, Ltd., Wellesley, MA (1997), xviii+572 pp

\bibitem{P} L. Paquet \emph{"Probl'emes mixtes pour le syst'eme de Maxwell"}, Annales Facult. des Sciences Toulouse (1982), Volume IV, 103-141

\bibitem{Re1} K. Reidemeister \emph{"Die Klassifikation der Linsenr\"{a}ume"}, Abhandl. Math. Sem. Hamburg 11 (1935), 102-109

\bibitem{Re2} K. Reidemeister \emph{"\"{U}berdeckungen von Komplexen"}, J. reine angew. Math. 173 (1935), 164-173 

\bibitem{Rh} G. de Rham \emph{"Complexes a automorphismes et homeomorphie differentiable"}, Ann. Inst. Fourier 2 (1950), 51-67

\bibitem{RaySin:RTA} D.B. Ray and I.M. Singer \emph{"R-Torsion and the Laplacian on Riemannian manifolds"}, Adv. Math. 7 (1971), 145-210

\bibitem{S} M. Spreafico \emph{"Zeta function and regularized determinant on a disc and on a cone"}, J. Geom. Phys. 54 (2005), no. 3, 355--371

\bibitem{Spr:ZIF} M. Spreafico \emph{"Zeta invariants for Dirichlet series"},  Pacific J. Math.  224  (2006),  no. 1, 185-200

\bibitem{Ver:ATO} B. Vertman \emph{"Analytic torsion of a bounded generalized cone"}, Comm. Math. Phys. 290 (2009), no. 3, 813--860. 

\bibitem{Ver:ZDF} B. Vertman \emph{"Zeta determinants for regular-singular Laplace-type operators"}, J. Math. Phys. 50 (2009), no. 8, 23 pp.

\bibitem{V1} B. Vertman \emph{"The Metric Anomaly at the Regular Boundary of the Analytic Torsion of a Bounded Generalized Cone, 
I. Odd-Dimensional Generalized Cone"}, arXiv:1004.2067v1 [math.SP] (2010)

\bibitem{V2} B. Vertman \emph{"The Metric Anomaly at the Regular Boundary of the Analytic Torsion of a Bounded Generalized Cone, 
II. Even-Dimensional Generalized Cone"},  arXiv:1004.2069v1 [math.SP] (2010)

\bibitem{W2} J. Weidmann \emph{"Linear Operators in Hilbert spaces"}, Graduate Texts in Mathematics, 68. Springer-Verlag, New York-Berlin, 1980. xiii+402 pp

\bibitem{Vi} S. Vishik \emph{"Generalized Ray-Singer Conjecture I. A manifold with smooth boundary"}, Comm. Math. Phys. 167 (1995), 1-102

\end{thebibliography}

\end{document}